\documentclass[]{article}
\usepackage{proceed2e}
\usepackage{amsmath,amssymb,amsfonts,amsthm}
\usepackage{bm}
\usepackage[hidelinks]{hyperref}

\usepackage{graphicx,xcolor} 
\usepackage[margin=1in]{geometry}
\usepackage{float}
\usepackage{times}
\usepackage{booktabs}
\usepackage{multirow, graphicx, rotating}
\usepackage{natbib}

\theoremstyle{definition}

\newtheorem{theorem}{Theorem}[section]
\newtheorem{corollary}[theorem]{Corollary}

\newtheorem{proposition}[theorem]{Proposition}

\makeatletter
\g@addto@macro\normalsize{%
  \setlength{\abovedisplayskip}{3pt}%
  \setlength{\belowdisplayskip}{3pt}%
  \setlength{\abovedisplayshortskip}{2pt}%
  \setlength{\belowdisplayshortskip}{2pt}%
}
\makeatother

\usepackage{titlesec}

\titlespacing*{\section}{0pt}{10pt}{6pt}
\titlespacing*{\subsection}{0pt}{8pt}{4pt}

\title{A probabilistic score estimation using the Bismut--Elworthy--Li formula}

\author{
  Jiatong Sun \thanks{~~Department of Mathematics and Statistics, University of Massachusetts Amherst}   \quad Long Nguyen \footnotemark[1] \quad
  Ning Jiang \footnotemark[1] \quad Weiqi Chu \footnotemark[1]\quad {\bf Yao Li} \thanks{~~Corresponding author} \, \footnotemark[1]
}

\begin{document}

\maketitle

\begin{abstract}
Score-learning methods generally lack pointwise supervision, which affects the training accuracy, particularly in regions where SDE trajectories are rare. To address this challenge, we introduce a novel method called \emph{BEL anchors}, which uses Monte Carlo score estimates at user-selected space-time query points. By combining a Feynman--Kac representation of the Fokker--Planck equation with the Bismut--Elworthy--Li formula, our estimator requires neither closed-form transition densities, spatial discretization, nor derivatives of the initial density. We derive its asymptotic variance and use the anchors as sparse supervised targets for existing score-learning objectives. Experiments show that BEL anchors improve implicit score matching in poorly sampled regions, prevent score-PINN collapse, and enable accurate sample regeneration in several numerical examples.
\end{abstract}

\section{Introduction} 
The score function
\[
\mathbf{s}(t,\mathbf{x})
=
\nabla_{\mathbf{x}}\log p_t(\mathbf{x})
\]
is a central object in score-based generative modeling
\citep{sohl2015deep,song2019generative,ho2020denoising,
song2021score} and in emerging numerical methods for
high-dimensional Fokker--Planck equations
\citep{boffi2023probability,hu2025score}.
In generative modeling, the score determines the reverse-time
dynamics that transport a reference distribution to the data
distribution. In scientific applications, it also provides access
to probability flows, evolving densities, and related physical
quantities such as entropy ption. Accurately learning the score is therefore important
beyond unconditional sample generation.

Most score-learning methods are unsupervised: pointwise values of
the true score are unavailable during training. Denoising score
matching replaces the marginal score with the conditional score of
a prescribed noising process \citep{vincent2011connection}, whereas
implicit and sliced score matching (ISM and SSM) eliminate the unknown score
through integration by parts
\citep{hyvarinen2005estimation,song2020sliced}. These objectives are
powerful, but their training objective is weighted by the distribution of
the training trajectories. Consequently, regions poorly covered by the trajectories contribute
little to the training loss, even when pointwise accuracy there is
important for a downstream task. Such regions arise, for example,
in long-time Langevin sampling and rare transition-path problems
\citep{cao2026robustness,triplett2025diffusion}.

A complementary approach is to enforce the score Fokker--Planck
equation through a physics-informed loss
\citep{lai2023fp,hu2025score}. This avoids relying entirely on
trajectory coverage, but a small PDE residual need not provide
effective supervision. In particular, the zero score
field exactly satisfies the interior score-PDE residual. When the
initial condition is imposed only through a soft penalty, optimization
can therefore collapse toward this trivial solution instead of
propagating the initial information through the space-time domain.

In this work, we introduce the \emph{BEL anchor}, a pointwise Monte
Carlo score estimator designed to supplement existing score-learning
objectives. For a user-selected query $(t,\mathbf{x})$, we first
rewrite the Fokker--Planck equation as a backward Feynman--Kac
equation for an auxiliary diffusion started from $\mathbf{x}$
\citep{milstein2004transition}. We then apply the
Bismut--Elworthy--Li formula
\citep{bismut1984large,elworthy1994formulae} to represent both
$p_t(\mathbf{x})$ and $\nabla_{\mathbf{x}}p_t(\mathbf{x})$ as
expectations over the same simulated paths. Their ratio yields a
pointwise estimator of the score. The construction requires neither
a spatial mesh, a conditional transition score, nor derivatives of
the initial density.

The BEL formula is a classical tool for differentiating Markov
semigroups and has been used to estimate gradients of
value functions in stochastic control and Hamilton--Jacobi--Bellman
equations \citep{bismut1978approximation,bismut1984large, sanders2025numerical}.
More recently, Malliavin--Bismut identities have been used to express
marginal or conditional diffusion scores through conditional
expectations, leading to regression-based score or control estimators
\citep{pidstrigach2025conditioning, mirafzali2025malliavin, mirafzali2025malliavin2}.
Our construction differs in both the represented object and the
resulting computational procedure.
We first transform the forward Fokker--Planck equation into an adjoint
Feynman--Kac problem whose auxiliary diffusion is initialized at the
query point, and then apply the BEL formula to obtain the pointwise
ratio estimator at query. This query-centered representation supplies sparse supervised score
anchors at user-selected locations, rather than learning the entire
score through a trajectory-weighted conditional regression.

The BEL anchor is not intended to replace existing score learning methods such as ISM and score-PINN. Instead, it supplies supervised targets at selected
space-time collocation points, whose distribution need not follow
the forward trajectories. We derive the BEL representation for
time-dependent SDEs with multiplicative noise and characterize the
asymptotic covariance of the resulting ratio estimator. The analysis
separates the contribution of density overlap from the growth of the
Malliavin weight and motivates rapidly mixing divergence-free
dynamics in place of the Ornstein--Uhlenbeck process. Numerically, we
show that sparse BEL anchors prevent the collapse of score-PINN and
improve the performance of ISM in regions poorly covered by its
training trajectories. Finally, we demonstrate accurate sample
regeneration in a few numerical examples.

\section{Methodology}
\label{sec:method}

\subsection{Probabilistic representation of the score}
\label{sec:fk}

Let $\mathbf{X}_t$ be a $d$-dimensional stochastic process given by
\begin{equation}
\label{eq:sde}
    d\mathbf{X}_t = \mathbf{f}(t,\mathbf{X}_t)\,dt + \boldsymbol{\sigma}(t,\mathbf{X}_t)\,d\mathbf{W}_t, \quad \mathbf{X}_0\sim p_0,
\end{equation}
where $\mathbf{f}$ is the drift, $\boldsymbol{\sigma}$ is the diffusion coefficient, and $\mathbf{W}_t$ is a standard $m$-dimensional Brownian motion. We write $\mathbf{a}=\boldsymbol{\sigma}\boldsymbol{\sigma}^\top$ as the diffusion matrix. Throughout, we assume that $\mathbf{a}$ is uniformly nondegenerate. In addition, we assume that $\mathbf{f}$ and $\boldsymbol{\sigma}$ has sufficient regularity such that the density $p(t,\cdot)$ of $\mathbf{X}_t$ is smooth and strictly positive for every $t>0$, so the score function $\mathbf{s}(t,\mathbf{x})=\nabla_{\mathbf{x}}\log p(t,\mathbf{x})$ is well defined. 

The density solves the Fokker--Planck equation 
\begin{equation}
    \label{eq:FPE}
    \partial_t p = -\nabla_{\mathbf{x}}\cdot(\mathbf{f}p) + \frac{1}{2}\sum_{i,j=1}^d\partial^2_{x_ix_j}(a_{ij}p)
\end{equation}
with $p(0,\cdot)=p_0$. Expanding the derivatives and using the symmetry of $\mathbf{a}$, we rewrite equation \eqref{eq:FPE} in the non-divergence form as
\begin{equation}
\label{eq:fpe-nondiv}
    \partial_t p = \boldsymbol{\beta}\cdot\nabla_{\mathbf{x}}p + \tfrac{1}{2}\,\mathbf{a}:\nabla^2_{\mathbf{x}}p + c\,p,
\end{equation}
where $\mathbf{A}:\mathbf{B}=\sum_{i,j}A_{ij}B_{ij}$, and the auxiliary drift $\boldsymbol{\beta}$ and the potential $c$ are given by
\begin{equation}
\label{eq:auxiliary-drift}
\begin{aligned}
    &\beta_i = -f_i + \sum_{j=1}^d\partial_{x_j}a_{ij}, \quad i=1,\ldots,d, \\
    &c = -\sum_{i=1}^d\partial_{x_i}f_i + \frac{1}{2}\sum_{i,j=1}^d\partial^2_{x_ix_j}a_{ij}.
\end{aligned}
\end{equation}

After a reversal of time, \eqref{eq:fpe-nondiv} becomes a backward Kolmogorov equation, which admits a Feynman--Kac probabilistic representation. Fix $t\in(0,T]$ and $\mathbf{x}\in\mathbb{R}^d$, and let $\{\mathbf{Y}_s^{t,\mathbf{x}}\}_{0\le s\le t}$ solve the reverse auxiliary SDE
\begin{equation}
\label{eq:reverse-auxiliary-sde}
    d\mathbf{Y}_s^{t,\mathbf{x}} = \boldsymbol{\beta}\bigl(t-s,\mathbf{Y}_s^{t,\mathbf{x}}\bigr)\,ds + \boldsymbol{\sigma}\bigl(t-s,\mathbf{Y}_s^{t,\mathbf{x}}\bigr)\,d\mathbf{W}_s
\end{equation}
where $\mathbf{Y}_0^{t,\mathbf{x}} =\mathbf{x}$ and $\mathbf{W}_s$ is a standard $m$-dimensional Brownian motion. We define the Feynman--Kac payoff $G_t^{t,\mathbf{x}}$ as
\begin{equation}
\label{eq:fk-payoff-general}
    G_t^{t,\mathbf{x}} = p_0\bigl(\mathbf{Y}_t^{t,\mathbf{x}}\bigr)\exp\left(\int_0^t c\bigl(t-s,\mathbf{Y}_s^{t,\mathbf{x}}\bigr)\,ds\right).
\end{equation}

The following Proposition, originally proposed in \cite{milstein2004transition}, provides the stochastic representation of Feynman--Kac solutions. See Appendix \ref{sec: FP_FK} for the proof. 

\begin{proposition}
    \label{prop:reverse-fk-density}
For every $t\in(0,T]$ and $\mathbf{x}\in\mathbb{R}^d$, the density of \eqref{eq:sde} satisfies
\begin{equation}
\label{eq:rep-fpe}
    p(t,\mathbf{x}) = \mathbb{E}\bigl[G_t^{t,\mathbf{x}}\bigr].
\end{equation}
\end{proposition}
Proposition~\ref{prop:reverse-fk-density} provides a Monte Carlo estimator for $p(t,\mathbf{x})$, but the score also requires $\nabla_{\mathbf{x}}p(t,\mathbf{x})$. It is impractical to differentiate the noisy Monte Carlo estimates by finite differences. Instead, we provide a probabilistic representation of the probability gradient $\nabla_{\mathbf{x}}p(t,\mathbf{x})$ using the Bismut--Elworthy--Li (BEL) formula, which transfers the derivative with respect to the starting point $\mathbf{x}$ to a stochastic weight.

Let $\mathbf{J}_s^{t,\mathbf{x}}=\partial_{\mathbf{x}}\mathbf{Y}_s^{t,\mathbf{x}}$ be the tangent process\footnote{See Appendix~\ref{sec:tangent-process} for details on the tangent process.} of \eqref{eq:reverse-auxiliary-sde}, where $\partial_{\mathbf{x}}$ means taking the partial derivative with respect to the initial value $\mathbf{x}$.  When no ambiguity arises, we suppress the superscripts and write $\mathbf{Y}_s$ and $\mathbf{J}_s$. Along each trajectory, we write $\boldsymbol{\sigma}_s^t=\boldsymbol{\sigma}(t-s,\mathbf{Y}_s)$, $\mathbf{a}_s^t=\mathbf{a}(t-s,\mathbf{Y}_s)$, and $\boldsymbol{\beta}_s^t=\boldsymbol{\beta}(t-s,\mathbf{Y}_s)$. The tangent process satisfies
\begin{equation}
\label{eq:reverse-tangent-sde}
\begin{aligned}
    d\mathbf{J}_s &= \nabla_{\mathbf{x}}\boldsymbol{\beta}^t_s\,\mathbf{J}_s\,ds + \sum_{k=1}^m\nabla_{\mathbf{x}}(\boldsymbol{\sigma}_{\cdot k})^t_s\,\mathbf{J}_s\,dW_s^k, \\
    \mathbf{J}_0 &= \mathbf{I}_d,
\end{aligned}
\end{equation}
where $\boldsymbol{\sigma}_{\cdot k}$ denotes the $k$-th column of $\boldsymbol{\sigma}$. 
\begin{theorem}[BEL representation of the density gradient]
\label{thm:BEL}
Suppose the assumptions in Appendix~\ref{app:assumptions} hold; in particular, $\boldsymbol{\beta}$, $\boldsymbol{\sigma}$, and $c$ are continuously differentiable in space with sufficient boundedness and integrability. Then, for every $t\in(0,T]$ and $\mathbf{x}\in\mathbb{R}^d$,
\begin{equation*}
\label{eq:general-density-gradient}
    \nabla_{\mathbf{x}}p(t,\mathbf{x}) = \mathbb{E}\bigl[G_t^{t,\mathbf{x}}\,\mathbf{H}_t^{t,\mathbf{x}}\bigr],
\end{equation*}
where the BEL weight $\mathbf{H}_t^{t,\mathbf{x}}$ is given by
\begin{equation}
\label{eq:general-bel-weight}
\begin{aligned}
    \mathbf{H}_t^{t,\mathbf{x}} = &\frac{1}{t}\int_0^t\mathbf{J}_s^\top(\mathbf{a}_s^t)^{-1}\boldsymbol{\sigma}_s^t\,d\mathbf{W}_s \\
    &+ \int_0^t\frac{t-s}{t}\,\mathbf{J}_s^\top\nabla_{\mathbf{x}}c(t-s,\mathbf{Y}_s)\,ds.
\end{aligned}
\end{equation}
\end{theorem}
We show the proof of Theorem \ref{thm:BEL} in Appendix \ref{sec:proof_thm_BEL}.

Proposition~\ref{prop:reverse-fk-density} and Theorem~\ref{thm:BEL} immediately give the following representation of the score.
\begin{corollary}[Probabilistic representation of the score]
\label{cor:score-representation}
Under the assumptions of Theorem~\ref{thm:BEL}, the score function is given by
\begin{equation}
\label{eq:score-representation}
    \mathbf{s}(t,\mathbf{x}) = \frac{\mathbb{E}\bigl[G_t^{t,\mathbf{x}}\,\mathbf{H}_t^{t,\mathbf{x}}\bigr]}{\mathbb{E}\bigl[G_t^{t,\mathbf{x}}\bigr]}.
\end{equation}
\end{corollary}

The representation \eqref{eq:score-representation} yields a Monte Carlo estimator of the score function that requires neither spatial differentiation nor spatial discretization. To evaluate it at a point $(t,\mathbf{x})$, we simulate $N$ independent sample paths of the reverse auxiliary SDE \eqref{eq:reverse-auxiliary-sde} and the tangent process \eqref{eq:reverse-tangent-sde}. On each sample path, both processes are driven by the \emph{same} realization of $\mathbf{W}_s$, and we compute $G_t^{t,\mathbf{x}}$ and $\mathbf{H}_t^{t,\mathbf{x}}$ from that path using \eqref{eq:fk-payoff-general} and \eqref{eq:general-bel-weight}. 

\subsection{Asymptotic variance of the score estimator}
\label{sec:score-variance}

Fix $(t,\mathbf{x})$ and abbreviate $G=G_t^{t,\mathbf{x}}$ and $\mathbf{H}=\mathbf{H}_t^{t,\mathbf{x}}$. Let $\{(G_n,\mathbf{H}_n)\}_{n=1}^N$ be independent copies of $(G,\mathbf{H})$, and define the Monte Carlo estimators of the density, its gradient, and the score by
\begin{equation}
\label{eq:mc-score-estimator}
\begin{aligned}
    \widehat{p}_N &= \frac{1}{N}\sum_{n=1}^N G_n, \quad 
    \widehat{\mathbf{g}}_N = \frac{1}{N}\sum_{n=1}^N G_n\mathbf{H}_n, \\
    \widehat{\mathbf{s}}_N &= \frac{\widehat{\mathbf{g}}_N}{\widehat{p}_N} = \frac{\sum_{n=1}^N G_n\mathbf{H}_n}{\sum_{n=1}^N G_n}.
\end{aligned}
\end{equation}
The asymptotic behavior of $\widehat{\mathbf{s}}_N$ follows from the strong law of large numbers and the delta method in the following theorem, and we show its proof in Appendix \ref{sec:proof_variance}. 
\begin{theorem}[Asymptotic variance of the score estimator]
\label{prop:score-asymptotic-variance}
Assume that $p(t,\mathbf{x})=\mathbb{E}[G]>0$, $\mathbb{E}[G^2]<\infty$, and $\mathbb{E}\bigl[G^2\|\mathbf{H}\|^2\bigr]<\infty$. Then $\widehat{\mathbf{s}}_N\to\mathbf{s}(t,\mathbf{x})$ almost surely as $N\to\infty$, and
\begin{equation}
\label{eq:score-clt}
\begin{aligned}
  &  \sqrt{N}\bigl(\widehat{\mathbf{s}}_N-\mathbf{s}(t,\mathbf{x})\bigr)\xrightarrow{\mathrm{d}}\mathcal{N}\bigl(\mathbf{0},\boldsymbol{\Sigma}(t,\mathbf{x})\bigr),
\end{aligned}\end{equation}
where
\begin{equation*}
\boldsymbol{\Sigma}(t,\mathbf{x}) = \frac{\mathbb{E}\bigl[G^2\bigl(\mathbf{H}-\mathbf{s}(t,\mathbf{x})\bigr)\bigl(\mathbf{H}-\mathbf{s}(t,\mathbf{x})\bigr)^\top\bigr]}{\mathbb{E}[G]^2}.
\end{equation*}
\end{theorem}
The asymptotic covariance has a useful factorization. Let $\mathbb{P}$ denote the probability measure under which $\mathbf{W}_s$ is a Brownian motion, and define the overlap factor $R_2$ and the tilted measure $\mathbb{Q}^{(2)}$ by
\begin{equation}
\label{eq:overlap-factor}
    R_2(t,\mathbf{x}) = \frac{\mathbb{E}[G^2]}{\mathbb{E}[G]^2}, \quad \frac{d\mathbb{Q}^{(2)}}{d\mathbb{P}} = \frac{G^2}{\mathbb{E}[G^2]}.
\end{equation}
Then
\begin{equation*}
    \boldsymbol{\Sigma}(t,\mathbf{x}) = R_2(t,\mathbf{x})\mathbb{E}_{\mathbb{Q}^{(2)}}\bigl[\bigl(\mathbf{H}-\mathbf{s}(t,\mathbf{x})\bigr)\bigl(\mathbf{H}-\mathbf{s}(t,\mathbf{x})\bigr)^\top\bigr].
\end{equation*}
The overlap factor $R_2$ measures the variability of the Feynman--Kac payoff, and the second term $\mathbb{E}_{\mathbb{Q}^{(2)}}$ measures the variability of the BEL weight under the measure tilted by $G^2$. 

In practice, the asymptotic covariance is estimated from the same samples by
\begin{equation*}
\label{eq:empirical-score-asymptotic-covariance}
    \widehat{\boldsymbol{\Sigma}}_N = \frac{1}{\widehat{p}_N^{\,2}}\,\frac{1}{N}\sum_{n=1}^N G_n^2\bigl(\mathbf{H}_n-\widehat{\mathbf{s}}_N\bigr)\bigl(\mathbf{H}_n-\widehat{\mathbf{s}}_N\bigr)^\top,
\end{equation*}
so that the covariance of $\widehat{\mathbf{s}}_N$ is approximately $\widehat{\boldsymbol{\Sigma}}_N/N$. 

\subsection{Score-based learning with BEL anchors}
\label{sec:bel-anchors}

The main application of the BEL score estimator is to provide a supervised component for existing score-based learning methods. Methods such as implicit score matching (ISM) \citep{hyvarinen2005estimation} and score-PINN \citep{hu2025score} constrain the network only through a variational objective or a PDE residual, together with an initial condition in the latter case. Neither observes pointwise values of the score for $t>0$. The BEL estimates \eqref{eq:mc-score-estimator} supply this missing information, and we show in Section~\ref{sec:ablation_study} that they can substantially improve training.

Let $\mathbf{s}_{\boldsymbol{\theta}}$ be a neural network approximation of the score with parameters $\boldsymbol{\theta}$, and let $\mathcal{L}_{\mathrm{base}}(\boldsymbol{\theta})$ be the loss of the underlying score-learning method. Given anchor points $\{(t_i,\mathbf{x}_i)\}_{i=1}^{N_{\mathrm{sup}}}$ and the corresponding BEL estimates $\widehat{\mathbf{s}}_i$ computed by \eqref{eq:mc-score-estimator}, we train with the combined loss
\begin{equation}
\label{eq:loss_combine}
    \mathcal{L}(\boldsymbol{\theta}) =  \mathcal{L}_{\mathrm{base}}(\boldsymbol{\theta}) + \frac{\lambda_{\mathrm{BEL}}}{N_{\mathrm{sup}}}\sum_{i=1}^{N_{\mathrm{sup}}}\bigl\|\mathbf{s}_{\boldsymbol{\theta}}(t_i,\mathbf{x}_i)-\widehat{\mathbf{s}}_i\bigr\|^2,
\end{equation}
where $\lambda_{\mathrm{BEL}}>0$ is a weighting parameter. In the numerical examples, the choice of $\mathcal{L}_{\mathrm{base}}$ depends on how the initial density is specified.

\paragraph{Initial density given by samples.}
When $p_0$ is available only through samples, we use ISM as the base method. Samples of $p_t$ are obtained by simulating \eqref{eq:sde} forward from the initial samples. In the payoff \eqref{eq:fk-payoff-general}, the unknown $p_0$ is replaced by a kernel density estimate $p_0^{\mathrm{KDE}}$ built from the same initial samples,
\begin{equation*}
\label{eq:fk-payoff-kde}
\begin{aligned}
    G_t^{t,\mathbf{x}} \approx{}& p_0^{\mathrm{KDE}}\bigl(\mathbf{Y}_t^{t,\mathbf{x}}\bigr) \exp\left(\int_0^t c\bigl(t-s,\mathbf{Y}_s^{t,\mathbf{x}}\bigr)\,ds\right).
\end{aligned}
\end{equation*}
Since the BEL weight \eqref{eq:general-bel-weight} does not involve $\nabla p_0$, no derivative of $p_0^{\mathrm{KDE}}$ is required. With this replacement, the BEL estimator targets the score of \eqref{eq:sde} initialized at $p_0^{\mathrm{KDE}}$ rather than $p_0$; the resulting discrepancy is controlled by the kernel bandwidth. The base loss is
\begin{equation*}
\label{eq:ism-loss}
\begin{aligned}
    \mathcal{L}_{\mathrm{ISM}}(\boldsymbol{\theta}) = \mathbb{E}_{t}\,\mathbb{E}_{\mathbf{x}\sim p_t}\Bigl[&\tfrac{1}{2}\bigl\|\mathbf{s}_{\boldsymbol{\theta}}(t,\mathbf{x})\bigr\|^2 +\nabla_{\mathbf{x}}\cdot\mathbf{s}_{\boldsymbol{\theta}}(t,\mathbf{x})\Bigr],
\end{aligned}
\end{equation*}
where $t\sim\mathcal{U}(0,T]$.

\paragraph{Initial density given explicitly.}
When $p_0$ is known in closed form but difficult to sample, we use score-PINN as the base method. We devide \eqref{eq:fpe-nondiv} by $p$ and take the gradient, showing that the score satisfies $\mathcal{L}^{\mathrm{PDE}}\mathbf{s}=\mathbf{0}$, where
\begin{equation*}
\label{eq:score-pde-operator}
\begin{aligned}
    \mathcal{L}^{\mathrm{PDE}}\mathbf{s} = \partial_t\mathbf{s} - \nabla_{\mathbf{x}}\Bigl[&c + \boldsymbol{\beta}\cdot\mathbf{s} + \tfrac{1}{2}\,\mathbf{a}:\bigl(\nabla_{\mathbf{x}}\mathbf{s}+\mathbf{s}\mathbf{s}^\top\bigr)\Bigr].
\end{aligned}
\end{equation*}
With collocation points drawn from a distribution $\mu$ on $(0,T]\times\mathbb{R}^d$ and initial points drawn from a distribution $\nu$ on $\mathbb{R}^d$, the base loss is
\begin{equation*}
\label{eq:score-PINN}
\begin{aligned}
    \mathcal{L}_{\mathrm{PINN}}(\boldsymbol{\theta}) = & \mathbb{E}_{(t,\mathbf{x})\sim\mu}\Bigl[\bigl\|\mathcal{L}^{\mathrm{PDE}}\mathbf{s}_{\boldsymbol{\theta}}(t,\mathbf{x})\bigr\|^2\Bigr] \\
    &+ \lambda_{\mathrm{IC}}\,\mathbb{E}_{\mathbf{x}\sim\nu}\Bigl[\bigl\|\mathbf{s}_{\boldsymbol{\theta}}(0,\mathbf{x})-\nabla_{\mathbf{x}}\log p_0(\mathbf{x})\bigr\|^2\Bigr],
\end{aligned}
\end{equation*}
where $\lambda_{\mathrm{IC}}>0$ weights the initial condition.

\section{Variance analysis and implementation}
\subsection{Variance of the Ornstein--Uhlenbeck drift} 
We first examine the overlap factor $R_2(t,\mathbf{x})$ for an
Ornstein--Uhlenbeck process. Consider the isotropic SDE
\begin{equation}
\label{eq:isotropic-ou}
    d\mathbf{X}_s
    =
    -\lambda\mathbf{X}_s\,ds
    +
    \sigma\,d\mathbf{B}_s,
    \quad
    \lambda>0,
\end{equation}
with Gaussian initial density $p_0 \sim N(0, \tau_0 \mathrm{I})$. For \eqref{eq:isotropic-ou}, the coefficients in
\eqref{eq:auxiliary-drift} are $
    \boldsymbol{\beta}(\mathbf{x})
    =
    \lambda\mathbf{x}, c=\lambda d $. Consequently, the auxiliary reverse process is 
\begin{equation}
\label{eq:reverse-ou}
    d\mathbf{Y}_s
    =
    \lambda\mathbf{Y}_s\,ds
    +
    \sigma\,d\mathbf{W}_s,
    \quad
    \mathbf{Y}_0=\mathbf{x}.
\end{equation}
The overlap factor can be explicitly calculated for this example. 

\begin{proposition}
\label{prop:ou-overlap-growth}
Let $\mathbf{Y}$ solve \eqref{eq:reverse-ou}, and let $p_0 \sim N(0, \tau_0 \mathrm{Id})$. Define $
    \nu_t
    =
    \frac{\sigma^2}{2\lambda}
    \left(e^{2\lambda t}-1\right),
    \quad
    r_t=\frac{\nu_t}{\tau_0^2}$.
Then
\begin{equation}
\label{eq:ou-R2-exact}
\resizebox{0.881\columnwidth}{!}{$\displaystyle
    R_2(t,\mathbf{x})=
    \left(
        \frac{1+r_t}{\sqrt{1+2r_t}}
    \right)^d
e^{\|\mathbf{x}\|^2e^{2\lambda t}
        \left(
            \frac{1}{\tau_0^2+\nu_t}
            -
            \frac{1}{\tau_0^2+2\nu_t}
        \right)}.
$}
\end{equation}
\end{proposition}

The overlap factor $R_2$ in equation \eqref{eq:ou-R2-exact} is exponential in both time and dimension. See the discussion after the proof of Proposition \ref{prop:ou-overlap-growth} in Appendix \ref{sec:ou-variance}. This rapidly amplifies the variance. In high dimension, this overlap deterioration means the BEL score estimate $\hat{\mathbf{s}}$ is only useful for very short $t$. 

\subsection{Analysis of divergence-free flow on the torus}
\label{sec:divergence-free-drift}

We next consider a divergence-free drift on the torus,
\begin{equation*}\label{eq:div-drift}
    \nabla_{\mathbf{x}}\cdot\mathbf{f}(\mathbf{x})=0,
\end{equation*}
together with additive isotropic noise
$\boldsymbol{\sigma}=\sigma\mathbf{I}_d$. In this case,
$\mathbf{a}=\sigma^2\mathbf{I}_d$, and the coefficients of the auxiliary
process reduce to $
    \boldsymbol{\beta}(\mathbf{x})=-\mathbf{f}(\mathbf{x}),
    \quad
    c(\mathbf{x})=0 $. The auxiliary process therefore satisfies
\begin{equation*}
\label{eq:div-free-auxiliary-sde}
    d\mathbf{Y}_s
    =
    -\mathbf{f}(\mathbf{Y}_s)\,ds
    +
    \sigma\,d\mathbf{W}_s,
    \quad
    \mathbf{Y}_0=\mathbf{x}.
\end{equation*}
Because $c=0$, the Feynman--Kac payoff contains no exponential
potential:
\begin{equation*}
\label{eq:div-free-payoff}
    G_t=p_0(\mathbf{Y}_t).
\end{equation*}
Moreover, the $\nabla_{\mathbf{x}}c$ term in the Malliavin weight
vanishes. Hence
\begin{equation}
\label{eq:div-free-malliavin-weight}
    \mathbf{H}_t
    =
    \frac{1}{\sigma t}
    \int_0^t
    \mathbf{J}_s^\top\,d\mathbf{W}_s,
\end{equation}
where $\mathbf{J}$ solves an ODE $
    \dot{\mathbf{J}}_s
    =
    -\nabla_{\mathbf{x}}\mathbf{f}(\mathbf{Y}_s)
    \mathbf{J}_s,
    \quad
    \mathbf{J}_0=\mathbf{I}_d$.
For a direction $\mathbf{v}\in\mathbb{R}^d$, the corresponding scalar
weight is
\begin{equation*}
\label{eq:div-free-directional-weight}
    H_{t,\mathbf{v}}
    =
    \mathbf{H}_t\cdot\mathbf{v}
    =
    \frac{1}{\sigma t}
    \int_0^t
    \mathbf{J}_s\mathbf{v}\cdot d\mathbf{W}_s.
\end{equation*}

Let $q_t^{\mathbf{x}}(\mathbf{y})$ denote the probability density function of
$\mathbf{Y}_t$ conditional on $\mathbf{Y}_0=\mathbf{x}$. Since the drift is divergence-free and the diffusion is spatially
constant, the uniform probability distribution on $\mathbb{T}^d$ is
invariant. Let $
    u(\mathbf{y})
    = |\mathbb{T}^d|^{-1} $
denote its density. If $q_t^{\mathbf{x}}$ is close to $u$, we have
\begin{equation*}
    R_2(t,\mathbf{x})
    \approx
    |\mathbb{T}^d|
    \int_{\mathbb{T}^d}
    p_0(\mathbf{y})^2\,d\mathbf{y}
    =
    \exp\left(
        D_2(p_0\|u)
    \right),
\end{equation*}
where
\begin{equation*}
    D_2(p_0\|u)
    =
    \log
    \int_{\mathbb{T}^d}
    \frac{p_0(\mathbf{y})^2}{u(\mathbf{y})}\,d\mathbf{y}
\end{equation*}
is the order-two R\'enyi divergence of $p_0$ from the uniform
distribution.  Thus, once the auxiliary process is close to equilibrium, the overlap
factor is determined primarily by the concentration of the initial
density relative to the uniform distribution.

The variance of the directional Malliavin weight follows from It\^o
isometry:
\begin{equation*}
    \operatorname{Var}(H_{t,\mathbf{v}})
    =
    \mathbb{E}\left[H_{t,\mathbf{v}}^2\right]=
    \frac{1}{\sigma^2t^2}
    \int_0^t
    \mathbb{E}\left[
        \|\mathbf{J}_s\mathbf{v}\|^2
    \right]\,ds.
\end{equation*}
This term isolates the contribution
of the expansion of the Jacobian matrix to the variance of the Malliavin weight. For large $t$, its growth rate depends on the Lyapunov exponent of the tangent flow \eqref{eq:reverse-tangent-sde}. 

Therefore, the variance of $H_{t,\mathbf{v}}$ is determined by two competing factors. On one hand, chaotic divergence-free flows may produce exponential growth of the tangent
process. On the other hand, the same advection can also substantially accelerate the
mixing of the noisy dynamics. The performance of the BEL score estimator depends on the comparison of those two exponential rates. 

\subsection{Practical implementation}
\label{sec:practical-implementation}
{\bf Choice of drift.} As discussed in Section \ref{sec:divergence-free-drift}, the performance of a divergence-free drift depends both on its Lyapunov exponent and its speed of convergence to equilibrium. The choice is not unique. In our simulations, we use a divergence-free vector field that is generalized from the Arnold–Beltrami–Childress (ABC) flow and additive white noise. Thus, equation \eqref{eq:sde} becomes
\begin{equation}
    \label{eq:ABC}
    \mathrm{d} x_i = [\sin(X_{i-1}) + \cos(X_{i+1})] dt + \sigma dW^i_t
\end{equation}
with cyclic indexing: $x_0 = x_d$, $x_{d+1} = x_1$ for $i = 1, \cdots, d$. The chaotic behavior of the ABC flow can be rigorously proved for $d = 3$. When $d > 3$, although a rigorous proof cannot be given, our numerical simulations show that the mixing time $t_{\mathrm{mix}}$ is stable as the dimension increases. In addition, our numerical simulations confirm that the variance of the Malliavin weight $\mathrm{Var}(\mathbf{H}_t)$ is approximately $\propto t^{-2}\exp( C(d) t)$, where the exponent $C(d)$ is also stable across all dimensions tested. We refer to Appendix \ref{app:mixing-vs-expansion} for the numerical results and further discussion of the two competing factors. 


{\bf Choice of collocation points.} 
 Without loss of generality, assume that a divergence-free drift term is used. If the score function $\mathbf{s}(t, \mathbf{x})$ is trained for sample regeneration, then collocation points are more useful at small $t$, because when $t$ is large, the probability density function is close to the uniform distribution on a torus. In sample regeneration examples, the sampling density of the temporal coordinate is concentrated near $t = 0$. We let the probability density at $t = 0$ be $10$ times of that at $t = 3$. As discussed in Appendix \ref{app:bel_strategies_spatial}, a uniform spatial distribution of collocation points is sufficient for most applications. We remark that this is different from the case of the Fokker--Planck solver \cite{zhai2022deep}.


{\bf Balance of different loss terms.}
The BEL anchor brings a new supervised loss function in addition to the original $\mathcal{L}_{\mathrm{base}}$. Therefore, it is important to balance the loss terms during neural network training. The general rule is to make the contributions of the different loss terms comparable. This can be done by scanning parameters for a suitable constant weight, using alternating ADAM proposed in \cite{zhai2022deep}, or calculating dynamic weights based on the gradient of each loss term in the training history as proposed in \cite{li2023artificial}. We find that a suitable constant weight is often good enough, as the relative scales of the loss terms do not change much during training.

\section{Numerical experiments}
\subsection{Ablation study}\label{sec:ablation_study}









The first experiment examines how the supervised BEL anchor changes the performance of existing score-learning methods. We report here the effect of changing the number of BEL collocation points. A second investigation of their spatial distribution is included in Appendix \ref{app:bel_strategies_spatial}. 

{\bf Common setup.} All ablation studies conducted in this section use the three-dimensional ABC flow on the torus $\mathbb T^3=[-\pi,\pi)^3$ as the forward dynamics, with the initial distribution given by a four-component Gaussian mixture model (GMM4). Both methods employ the same neural network architecture and BEL anchor, while their respective unsupervised objectives differ. The combined training losses are given by $\mathcal{L}_{\mathrm{ISM}}(\boldsymbol{\theta})$ and $\mathcal{L}_{\mathrm{PINN}}(\boldsymbol{\theta})$ as defined in Section~\ref{sec:bel-anchors}. We refer to Appendix~\ref{app:bel_number_spatial} for the full simulation details. 

The training set of BEL anchors consists of $N_{\mathrm{BEL}}$ collocation points sampled uniformly on $\mathbb{T}^3$. For each spatial point $\mathbf{x}_i$, we independently sample 10 temporal coordinates from $[0,3]$ uniformly and independently, bringing the total number of BEL collocation points to $10N_{\mathrm{BEL}}$.


{\bf Numerical reference.} Since the exact score function is unavailable, we construct two method-specific numerical reference scores. Specifically, the ISM reference uses ISM with BEL anchors, 40,000 stochastic trajectories for ISM and 40,000 spatial BEL collocation points. The neural network is trained for 60,000 steps. Similarly, the score-PINN reference is obtained using 20,000 uniformly distributed spatial collocation points for the PINN residual and 20,000 spatial collocation points for BEL anchors. The score-PINN reference network is trained for 40,000 optimization steps. The resulting score approximations, denoted by $s^{\mathrm{ISM}}_{\mathrm{ref}}(x,t)$ and $s^{\mathrm{PINN}}_{\mathrm{ref}}(x,t)$, are used as the numerical reference solutions for the remaining ablation studies.

{\bf Experiment results.} To examine the effect of BEL anchors, we progressively reduce the number of BEL spatial collocation points from 20,000 to 0, using $N_{\mathrm{BEL}} = 0,2,10,20,50,150,300,625,
1250,2500,5000,10000$, and $20000$, respectively.
For each positive $N_{\mathrm{BEL}}$, the space-time collocation points are sampled in the same way as in training the reference solution. The ISM trajectory dataset, PINN collocation dataset, neural network architecture, and other training settings are held fixed across the ablation experiments. When $N_{\mathrm{BEL}}=0$, both methods become unsupervised.

\begin{figure}[h!]
    \centering
    \includegraphics[width=0.9\linewidth]{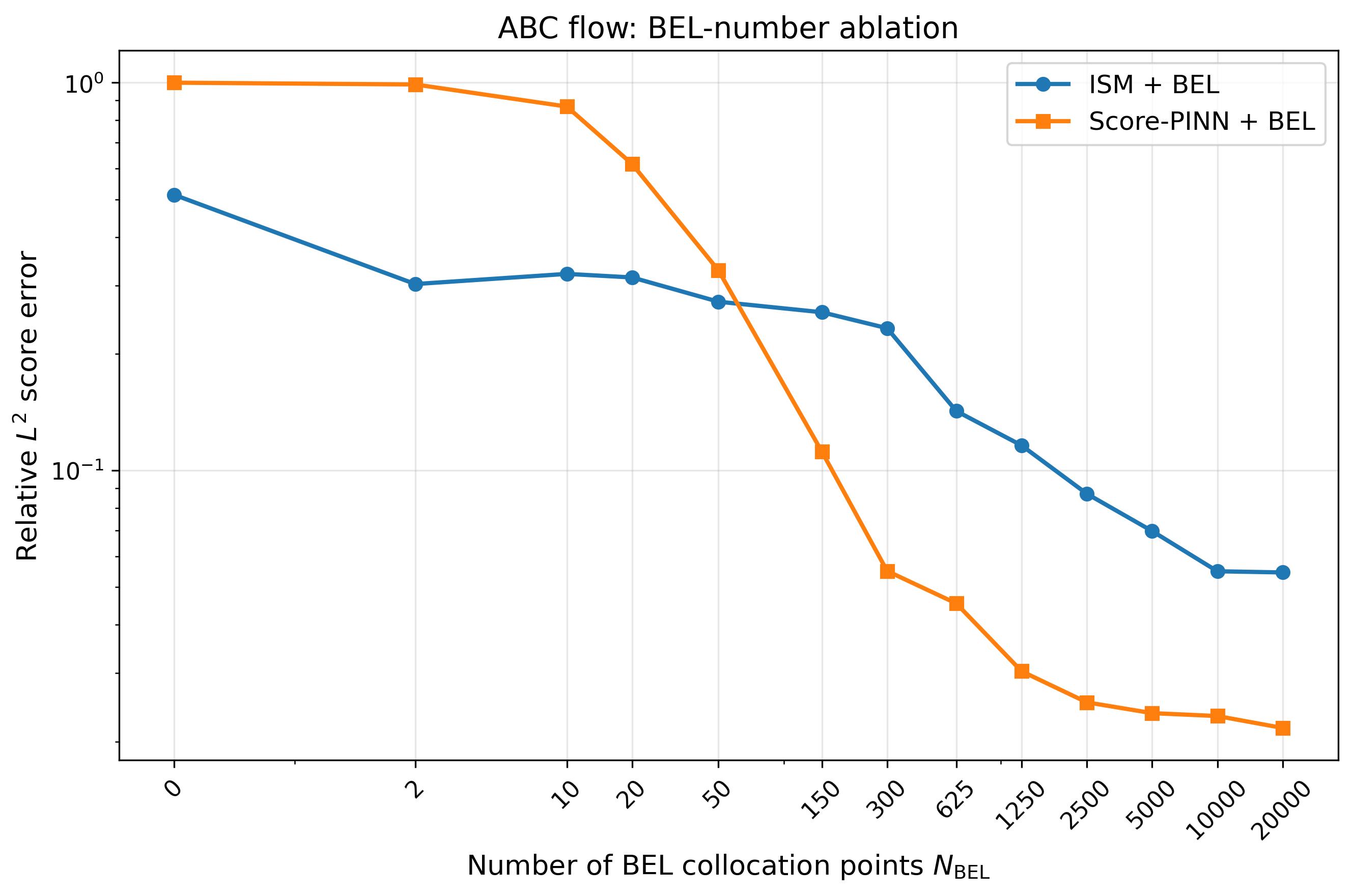}
    \caption{Log-log plots of the relative $L^2$ score errors of ISM with BEL anchor and score-PINN with BEL anchor as functions of the number of BEL spatial collocation points $N_{\mathrm{BEL}}$.}
    \label{fig:bel_number_error}
\end{figure}

The effect of BEL anchors is measured by the relative $L^2$ score error over the space-time domain
$$\small
E_{L^2}^{m}=\sqrt{\frac{\int_{0}^{1}\int_{\mathbb T^3}\|s_{\boldsymbol{\theta}}^{m}(x,t)-s_{\mathrm{ref}}^{m}(x,t)\|_2^2\,dx\,dt}{\int_{0}^{1}\int_{\mathbb T^3}\|s_{\mathrm{ref}}^{m}(x,t)\|_2^2\,dx\,dt}}\;,
$$
where $m\in\{\mathrm{ISM},\mathrm{PINN}\}$, evaluated against the method-specific numerical reference described above, using 20,000 spatial points sampled uniformly from $\mathbb T^3$ and 10 equally spaced time points $t_k=0.1k$, $k=1,\ldots,10$.
Figure~\ref{fig:bel_number_error} presents the resulting relative $L^2$ errors as the number of BEL collocation points decreases. A symmetric logarithmic scale is used for the horizontal axis to accommodate $N_{\mathrm{BEL}}=0$, while the vertical axis uses a logarithmic scale. The corresponding numerical values are reported in Table~\ref{tab:bel_number_error} in Appendix~\ref{app:bel_number_spatial}.

As shown in Figure~\ref{fig:bel_number_error} above and Table~\ref{tab:bel_number_error} in Appendix~\ref{app:bel_number_spatial}, the two methods exhibit markedly different sensitivities to the number of BEL collocation points. Score-PINN can tolerate a significantly smaller number of collocation points and produce a reasonable score estimate. However, once there are too few collocation points, score-PINN collapses to a trivial solution. This means standalone score-PINN cannot learn the score function, as pointed out in other literature \cite{hu2025score}. BEL anchors make this method applicable. Compared with score-PINN, ISM is more robust. As seen in Figure~\ref{fig:bel_heatmap_L2_error} in Appendix~\ref{app:bel_number_spatial}, ISM can still  approximate the score without any supervision in the regions with sufficient coverage of SDE trajectories. BEL anchors improve score quality in regions that random trajectories fail to cover during ISM training. In both cases, this ablation study demonstrates the effectiveness of BEL anchors in improving the existing score-learning methods. 

In addition, we examine the effect of the spatial distribution of BEL collocation points. We refer to Appendix~\ref{app:bel_strategies_spatial} for details.

\subsection{Variance scaling}
\label{sec:var-scaling}

\begin{figure}[H]
    \centering
    \includegraphics[width=\linewidth]{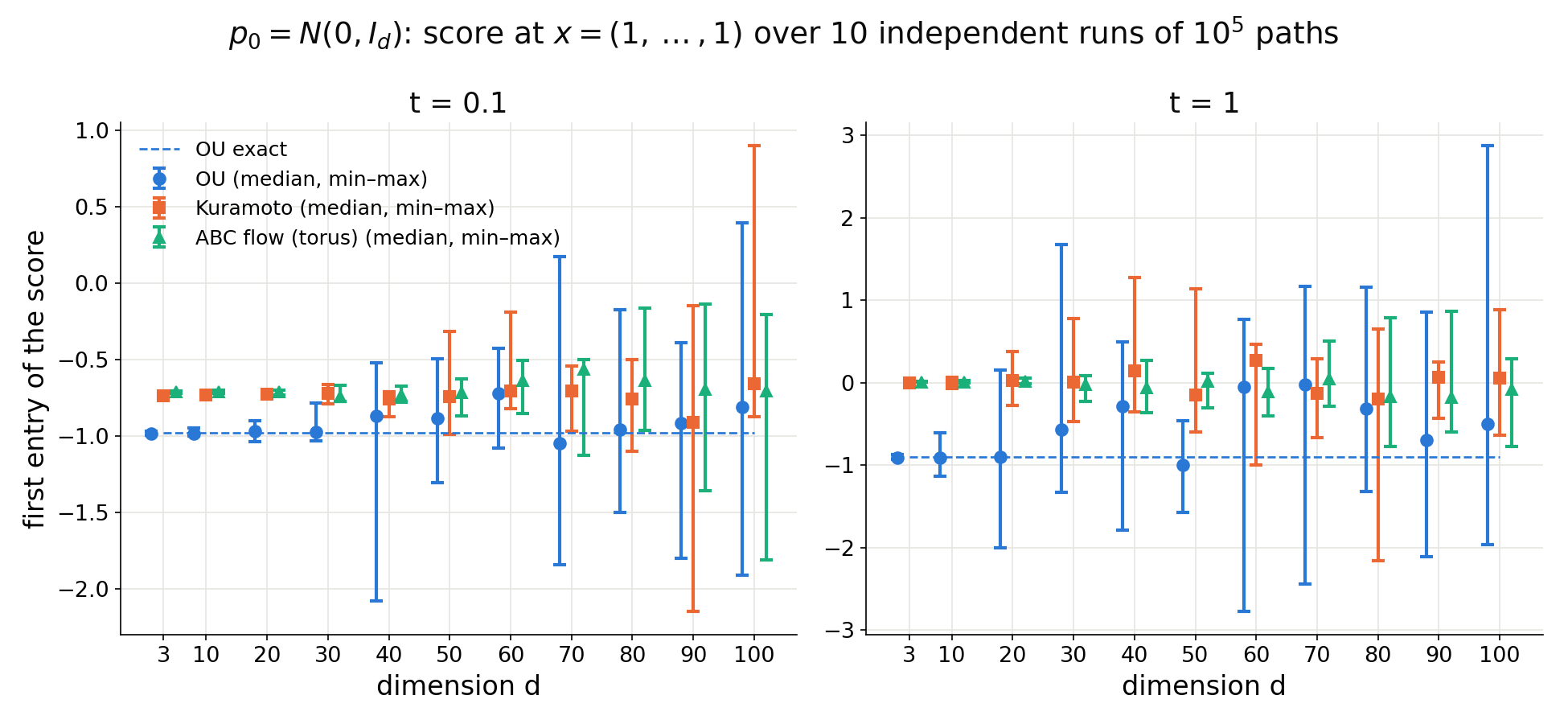}
    \includegraphics[width=\linewidth]{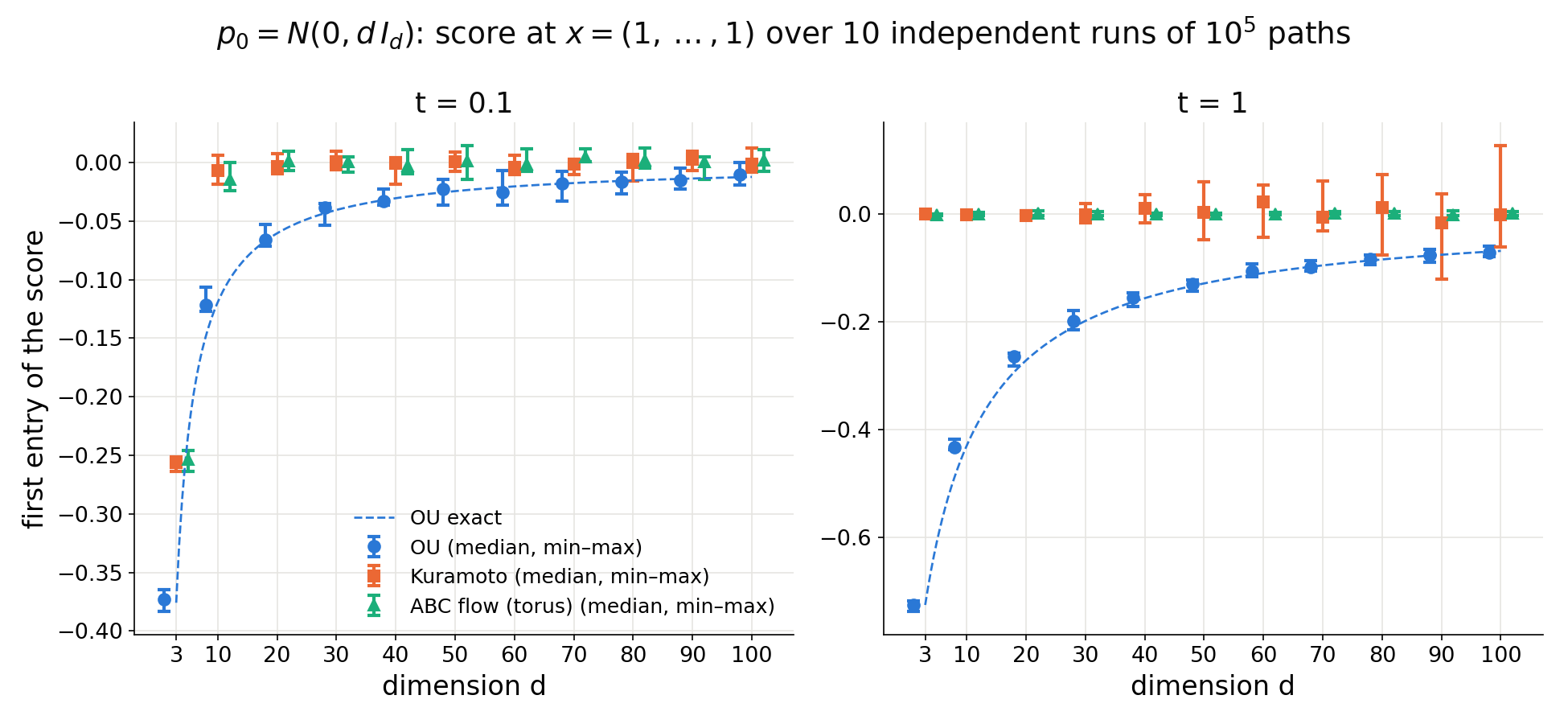}
    \caption{Median, minimum, and maximum of $10$ independent runs of score at $\mathbf{x} = (1, 1, \cdots, 1)$ at each dimension for three different underlying drifts.}
    \label{fig:var-1}
\end{figure}

In this experiment, we examine variance scaling as the dimension increases. We choose three underlying drifts: the OU process, the generalized ABC flow, and the Kuramoto process as described in Appendix~\ref{app:Variance_scaling}. To examine the effect of the overlap factor $R_2$, we consider two initial distributions $\mathcal{N}(0, I_d)$ and $\mathcal{N}(0, d I_d)$, since the $R_2$ factor for $\mathcal{N}(0, I_d)$ increases exponentially with dimension $d$. Then we compute the median, minimum and maximum of the first entry of the score at $\mathbf{x} = (1,1,...,1)$ over 10 independent runs using $10^5$ Monte Carlo samples. The results are shown in Figure \ref{fig:var-1}. 

Because a $d$-dimensional box is exponentially larger than a $d$-dimensional sphere, for the initial density $\mathcal{N}(0, I_d)$, the overlap factor $R_2$ grows exponentially in both $d$ and $t$. As seen in the top panel of Figure \ref{fig:var-1}, the estimated score values have high variance for all three drifts when the dimension is very high. The OU process has the worst scalability with the dimension, as discussed in Proposition \ref{prop:ou-overlap-growth}. The generalized ABC flow behaves better because it is restricted to a torus and has a more slowly expanding $\mathbf{H}_t$. As seen in Figure \ref{fig:var-4} in Appendix~\ref{app:Variance_scaling}, at $t = 1$, ABC flow can still produce about $100$ effective samples with $10^6$ trajectories in dimensions $30\sim 40$, which is sufficient for BEL anchor. Higher dimension many require multi-stage BEL estimator.

On the other hand, the initial density $\mathcal{N}(0, dI_d)$ expands together with the torus/box as its dimension $d$ increases. Hence its overlap factor $R_2$ does not grow exponentially with the dimension. As seen in the bottom panel of Figure \ref{fig:var-1}, the estimated score values have significantly lower variance. However, at $t = 1$, the score variance for the Kuramoto process is significantly higher because its $G_t$ contains an exponential term $\exp(\int_0^t c  \,\mathrm{d}s)$.

The variance scaling experiment confirms that the generalized ABC flow generally produces smaller variance, particularly when the initial distribution is concentrated. In very high-dimensional problems with very large overlap factor $R_2$, the BEL anchors must be computed in stages to further reduce the variance. We will explore this in our subsequent work.


\subsection{Generative performance}
We provide the following three examples to examine the performance of sample regeneration using our BEL-based score estimation.

\begin{figure}[]
    \centering
    \includegraphics[width=\linewidth]{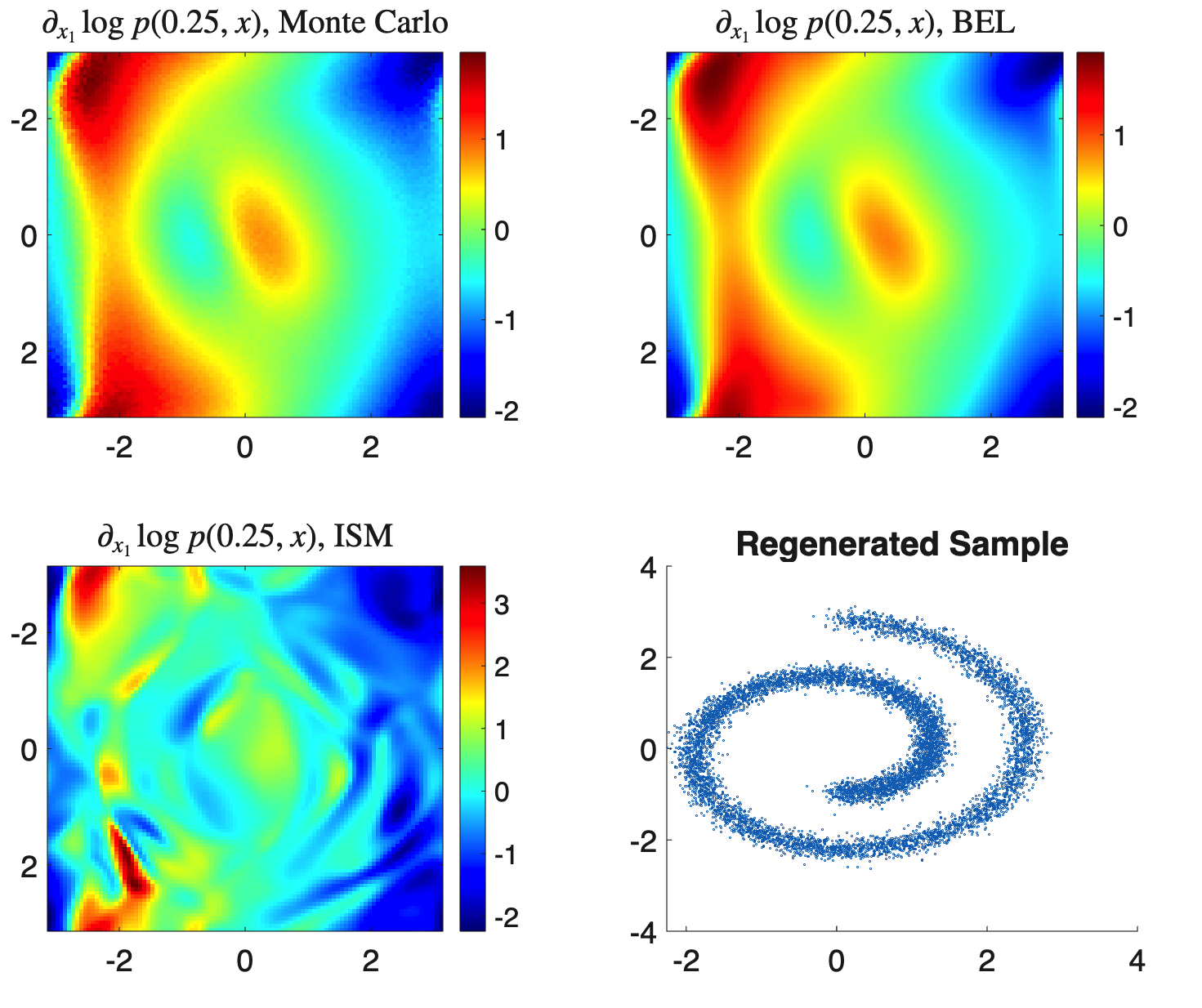}
    \caption{Top: Monte Carlo score estimate, learned score with BEL anchor. Bottom left: learned score only using ISM. Bottom right: Regenerated samples using learned score with BEL anchor.}
    \label{fig:swiss_roll}
\end{figure}

{\bf Swiss roll.} The first example is the Swiss roll dataset in 3D for a visual check. The initial dataset is concentrated on the manifold 
\begin{equation*}
\{ (x,y,z) \in \mathbb{T}^d\,|\, x = t \cos t, y = t\sin t, t \in [1.5 \pi, 4.5 \pi]\}
\end{equation*}
with isotropic noise $\delta = 0.1$. We use KDE to approximate the initial PDF $p_0$. To obtain a reference solution, we first use numerical simulation as in equation \eqref{eq:score-representation} to estimate the score function on three $(x,y,0)$ planes at $t = 0.25, 0.5, 1.0$ respectively. Then we draw $20000$ uniform samples from $\mathbb{T}^3$ as the training set. The BEL estimate at each point uses $10000$ independent trajectories to estimate the probability density function and its derivative at $10$ random times. Next, we use ISM as the base score-learning method to examine the effect of BEL anchors in learning the score function. The score functions learned with and without a BEL anchor are compared with the Monte Carlo estimates on the three aforementioned slices. The $L^2$ distances for the $3$ spatial components of the score function obtained with and without the BEL anchor are listed in Table \ref{tab:swiss_roll}. It is easy to see that the use of BEL anchors significantly improves the quality of the score function. As a visual confirmation, in Figure \ref{fig:swiss_roll}, we demonstrate a comparison of the learned $\mathbf{s}_1$ at $t = 0.25$ with or without the BEL anchor.


\begin{table}[h!]
    \centering
    \renewcommand{\arraystretch}{1.2}

    \resizebox{\columnwidth}{!}{%
        \begin{tabular}{ccc}
            \toprule
            \textbf{$t$} 
            & $\boldsymbol{(e_1,e_2,e_3)}$, \textbf{BEL}
            & $\boldsymbol{(e_1,e_2,e_3)}$, \textbf{ISM} \\
            \midrule
            0.25 & (0.0354, 0.0369, 0.0268)
                 & (0.4533, 0.3938, 0.2299) \\
            0.5  & (0.0153, 0.0206, 0.0154)
                 & (0.2157, 0.2241, 0.1872) \\
            1.0  & (0.0108, 0.0099, 0.0079)
                 & (0.1298, 0.1060, 0.0734) \\
            \bottomrule
        \end{tabular}%
    }

    \caption{$L^2$ distance between learned score function and Monte Carlo estimate, where $\mathbf{e}_i$ is the $L^2$ distance between learned $\mathbf{s}_i$ and Monte Carlo simulated $\mathbf{s}_i$.}
    \label{tab:swiss_roll}
\end{table}

As a visual confirmation of the sample regeneration result, we use the learned score function to regenerate the Swiss roll dataset, which is demonstrated in the right panel of Figure \ref{fig:swiss_roll}.

\begin{figure}[h!]
    \centering
    \includegraphics[width=\linewidth]{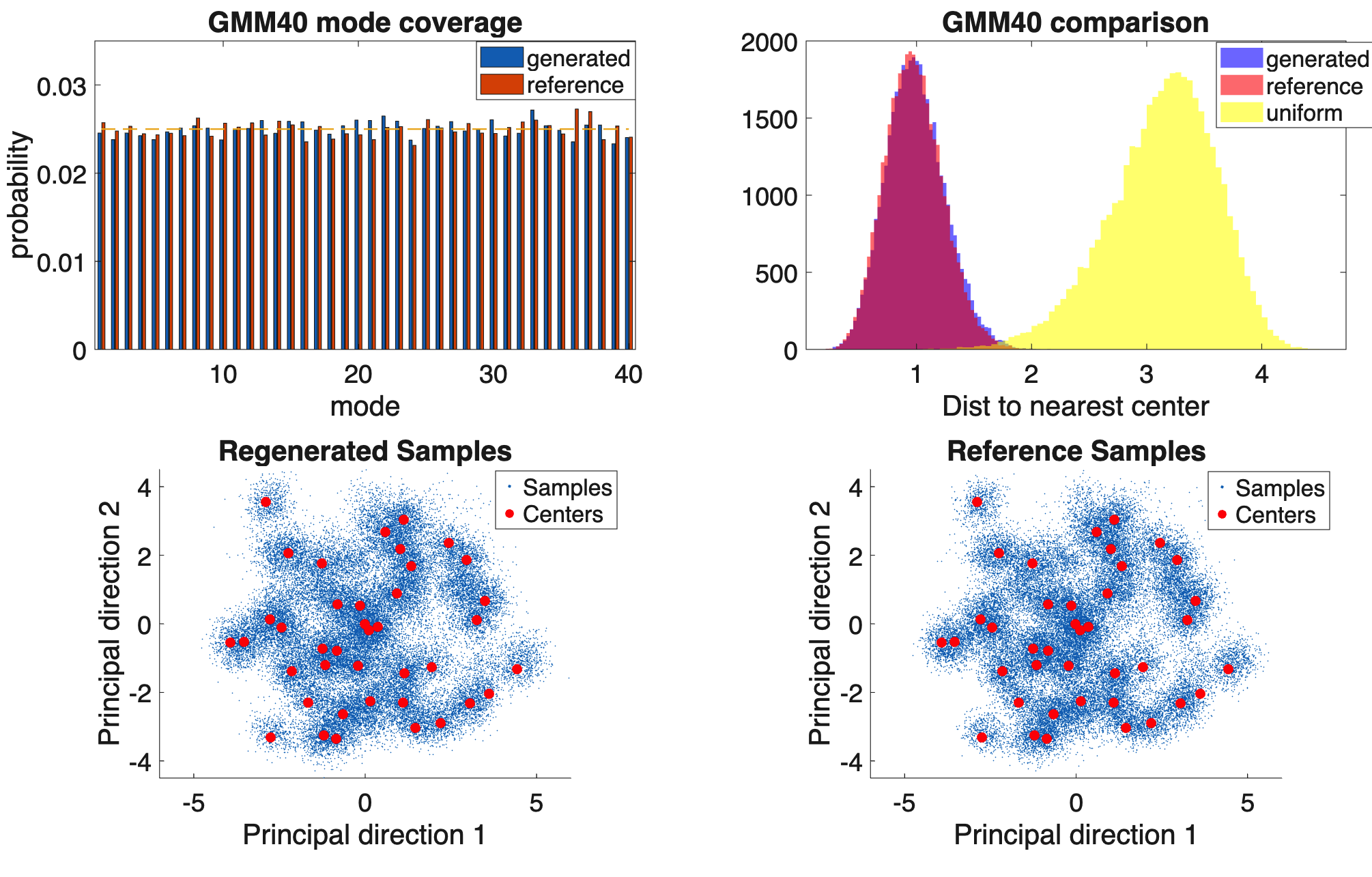}
    \caption{Top Left: Mode coverage of regenerated samples. Top Right: Distribution of the distance to the nearest center $\mathbf{y}_i$. Bottom: Regenerated samples and Reference samples, projected to the plane spanned by two principal directions of Gaussian centers.}
    \label{fig:gmm8}
\end{figure}

{\bf Gaussian mixture model.} The next example is the sample regeneration of a 40-mode Gaussian mixture in $\mathbb{T}^8$. We randomly sample $40$ centers, denoted by $\{ \mathbf{y}_1 ,\cdots, \mathbf{y}_{40}\}$ in $\mathbb{T}^8$. Let $p_0$ be the Gaussian mixture 
\begin{equation*}\label{eq:40D-gaussian-mixture}
p_0(\mathbf{x}) = \frac{1}{40}\sum_{i = 1}^{40} \mathcal{N}(\mathbf{x}; \mathbf{y}_i, 0.35 \mathrm{I}) \,,
\end{equation*} 
where $\mathcal{N}$ is the probability density function of a normal random variable wrapped on the torus. Since $p_0$ is given as an explicit probability density rather than samples, implicit score matching cannot be used here. Instead, we use BEL anchors together with score-PINN. As before, we sample $20000$ collocation points in $\mathbb{T}^8$ and use the BEL formula to estimate the probability density function and its derivatives at $10$ random times. The BEL anchor is used together with score-PINN using a combined loss function as in equation \eqref{eq:loss_combine}. The learned score function is then used to regenerate the 40-mode Gaussian mixture from the uniform distribution on $\mathbb{T}^8$. With the help of BEL anchors, the 40-mode Gaussian mixture is successfully regenerated. The regenerated Gaussian mixture distribution, denoted by $\hat{p}_0$, is visually indistinguishable from the original $p_0$ (bottom two panels of Figure \ref{fig:gmm8}). The Sinkhorn distance between $p_0$ and  $\hat{p}_0$ is $0.0158$ when measured using $2000$ points, smaller than the Sinkhorn distance reported in literature \cite{choi2026non} using other methods. As seen in the top two panels of Figure \ref{fig:gmm8}, $p_0$ and $\hat{p}_0$ are consistent in both mode coverage and distance from the nearest center $\mathbf{y}_i$.

\begin{figure}[h!]
    \centering
    \includegraphics[width=\linewidth]{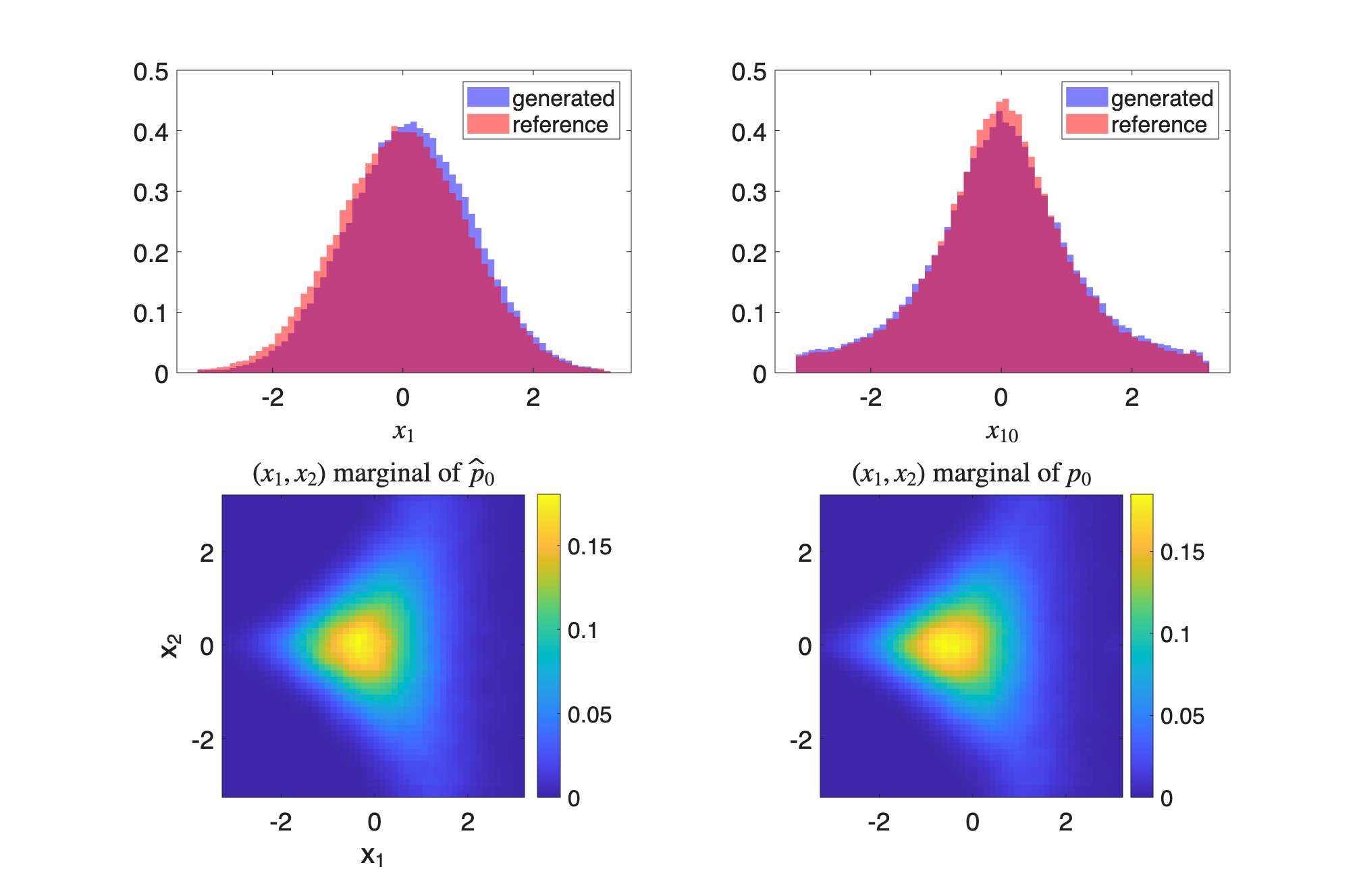}
    \caption{Top: Marginal distribution of $x_1$ and $x_{10}$ of generated samples and reference samples. Bottom: Heat map of $(x_1,x_2)$-marginal distributions of generated samples and reference samples.}
    \label{fig:funnel}
\end{figure}

{\bf Neal's funnel.} The last example is Neal's funnel in $10$D. Let the initial probability density function be
\begin{equation*}
\tilde{p}_0(x_1, \cdots, x_10) = \mathcal{N}(x_1;0,0)\mathcal{N}((x_2, \cdots, x_{10}); 0, e^{x_1})
\end{equation*}and $p_0$ be the wrapped probability density function of $\tilde{p}_0$ onto $\mathbb{T}^{10}$. The training set simulation and score-function training are similar to those in the GMM40 example. We still use score-PINN as the initial distribution is given as a probability density function rather than samples. The learned score function is then used to regenerate Neal's funnel in $10$D. The result is demonstrated in Figure~\ref{fig:funnel}. As seen in Figure \ref{fig:funnel}, the regenerated Neal's funnel, denoted by $\hat{p}_0$, is consistent with $p_0$. The Sinkhorn distance between $p_0$ and $\hat{p}_0$ is $0.0145$ when measured using $2000$ sample points. 

\section{Conclusion}

We introduced the BEL anchor, a pointwise Monte Carlo estimator that supplies supervised score information at user-selected space--time locations. The construction first represents the solution of the Fokker--Planck equation through an auxiliary Feynman--Kac diffusion initialized at the query point and then applies the Bismut--Elworthy--Li formula to its spatial derivative. This yields a ratio estimator for the score that does not require a closed-form transition density, a spatial mesh, or differentiation of the initial density. Consequently, the same construction applies when the initial distribution is specified either by an evaluable density or, after density estimation, by samples.

The BEL anchor is intended to complement rather than replace existing global score-learning methods. Our experiments demonstrate this complementarity. Implicit score matching (ISM) remains effective in well-sampled regions but can be substantially improved by anchors in regions that are poorly covered by forward trajectories. A score-PINN, by contrast, may collapse to an uninformative solution when trained only with its PDE residual and initial condition; a comparatively sparse set of BEL anchors removes this ambiguity and enables accurate learning. The resulting hybrid methods accurately regenerate samples in our numerical experiments. 

We also derived the asymptotic covariance of the score estimator and decomposed it into a density-overlap factor and the variability of the BEL weight under a tilted measure. This decomposition explains why apparently stable estimates can become unreliable when the effective overlap is small and why the choice of auxiliary drift matters. The analytical and numerical comparisons indicate that divergence-free mixing dynamics can offer substantially better variance behavior than an Ornstein--Uhlenbeck drift. Together, these results establish BEL anchors as a practical mechanism for adding targeted pointwise information to otherwise unsupervised or weakly constrained score-learning problems.

\subsubsection*{Acknowledgements}

WC was partially supported by NSF Awards DMS-2350325 and DMS-2514053. YL was partially supported by NSF award DMS-2510209.

\bibliographystyle{apalike}
\bibliography{references}

\clearpage

\appendix

\section{Probability Preliminary}\label{app:prelim}

\subsection{Standing assumptions}\label{app:assumptions}

We summarize the assumptions required by Section~\ref{sec:fk}. 

\paragraph{(A1) Existence of SDE solution.}
We assume that $\mathbf{f}$ and $\boldsymbol{\sigma}$ has sufficient regularity such that equation \eqref{eq:sde} has a unique strong solution with continuous paths and finite second moments at finite time.

\paragraph{(A2) Nondegeneracy.}
Suppose that $m\ge d$ and there exist constants $0<\lambda\le\Lambda<\infty$ with
\begin{equation*}
  \lambda\|\boldsymbol{\xi}\|^2
  \le \boldsymbol{\xi}^{\top}\mathbf{a}(t,\mathbf{x})\,\boldsymbol{\xi}
  \le \Lambda\|\boldsymbol{\xi}\|^2
\end{equation*}
for all $\boldsymbol{\xi}\in\mathbb{R}^d$ and all $(t,\mathbf{x})$. Let $\boldsymbol{\sigma}^{+}=\boldsymbol{\sigma}^{\top}\mathbf{a}^{-1}$ and assume $\sup_{\mathbf{x},t}\|\boldsymbol{\sigma}^{+}(t,\mathbf{x})\|<\infty$. Under (A1) and (A2), the law of $\mathbf{X}_t$ in \eqref{eq:sde} admits a density $p(t,\cdot)$ that is smooth and strictly positive on $\mathbb{R}^d$ for every $t>0$.

\paragraph{(A3) Initial condition.}
The state $\mathbf{X}_0$ is independent of $\mathbf{W}$ and its law admits a continuous, bounded, and strictly positive density $p_0$ with finite second moment. 

\paragraph{(A4) Integrability of the exponential weight.} 
With $c$ in \eqref{eq:auxiliary-drift}, we require $\mathbb{E}\big[e^{2c(\mathbf{X}_t)}\big]<\infty$ for every $t\in(0,T]$ and every starting point of interest. 

\subsection{The tangent process of a gradient flow}
\label{sec:tangent-process}
Fix $t\in(0,T]$. Under the assumptions of Appendix~\ref{app:assumptions}, the reverse auxiliary SDE \eqref{eq:reverse-auxiliary-sde} has a unique strong solution $\mathbf{Y}_s^{t,\mathbf{x}}$ for every initial point $\mathbf{x}\in\mathbb{R}^d$. Moreover, the map $\mathbf{x}\mapsto\mathbf{Y}_s^{t,\mathbf{x}}$ admits a version that is almost surely continuously differentiable, simultaneously for all $s\in[0,t]$ \citep{kunita1990stochastic}, \citep{protter2005stochastic}[Chapter~V]. The \emph{tangent process} is the Jacobian matrix of this map,
\begin{equation*}
\label{eq:tangent-def}
    \mathbf{J}_s^{t,\mathbf{x}} = \nabla_{\mathbf{x}}\mathbf{Y}_s^{t,\mathbf{x}},
\end{equation*}
with entries $(\mathbf{J}_s^{t,\mathbf{x}})_{ij}=\partial_{x_j}(Y_s^{t,\mathbf{x}})_i$ for $i,j=1,\ldots,d$. Its $j$-th column, $\mathbf{J}_s^{t,\mathbf{x}}\mathbf{e}_j$, measures the first-order response of the trajectory at time $s$ to a perturbation of the initial point in the direction $\mathbf{e}_j$.

We derive the governing SDE for the tangent process.
Writing \eqref{eq:reverse-auxiliary-sde} in integral form, we have
\begin{equation}
\label{eq:Y-integral}
\begin{aligned}
    \mathbf{Y}_s = &\mathbf{x}
    + \int_0^s\boldsymbol{\beta}\bigl(t-r,\mathbf{Y}_r\bigr)\,dr \\
    &+ \sum_{k=1}^m\int_0^s\boldsymbol{\sigma}_{\cdot k}\bigl(t-r,\mathbf{Y}_r\bigr)\,dW_r^k .
\end{aligned}
\end{equation}
Differentiating \eqref{eq:Y-integral} with respect to $x_j$, and interchanging differentiation with the Lebesgue and It\^o integrals as justified in \cite{kunita1990stochastic}, the chain rule gives
\begin{equation*}
\label{eq:J-column}
\begin{aligned}
    \partial_{x_j}\mathbf{Y}_s = &\mathbf{e}_j
    + \int_0^s\nabla_{\mathbf{x}}\boldsymbol{\beta}^t_r\,\partial_{x_j}\mathbf{Y}_r\,dr \\
    &+ \sum_{k=1}^m\int_0^s\nabla_{\mathbf{x}}(\boldsymbol{\sigma}_{\cdot k})^t_r\,\partial_{x_j}\mathbf{Y}_r\,dW_r^k ,
\end{aligned}
\end{equation*}
where $\nabla_{\mathbf{x}}\boldsymbol{\beta}^t_r$ and $\nabla_{\mathbf{x}}(\boldsymbol{\sigma}_{\cdot k})^t_r$ denote the $d\times d$ Jacobian matrices of $\boldsymbol{\beta}$ and $\boldsymbol{\sigma}_{\cdot k}$ in the spatial variable, evaluated at $(t-r,\mathbf{Y}_r)$. Collecting the columns $j=1,\ldots,d$ yields the matrix-valued linear SDE \eqref{eq:reverse-tangent-sde} with $\mathbf{J}_0=\mathbf{I}_d$. Since $\nabla_{\mathbf{x}}\boldsymbol{\beta}$ and $\nabla_{\mathbf{x}}\boldsymbol{\sigma}_{\cdot k}$ are bounded, the Burkholder--Davis--Gundy and Gronwall inequalities give the moment bounds
\begin{equation*}
\label{eq:J-moments}
    \mathbb{E}\Bigl[\sup_{0\le s\le t}\bigl|\mathbf{J}_s^{t,\mathbf{x}}\bigr|^p\Bigr] \le C_p, \qquad p\ge 1,
\end{equation*}
with $C_p$ independent of $\mathbf{x}$. These bounds ensure that the BEL weight \eqref{eq:general-bel-weight} has finite moments of all orders.

\subsection{Malliavin derivative and Skorokhod integral}
\label{sec:malliavin}

We briefly recall the elements of Malliavin calculus needed in our
analysis. Let
$(\Omega,\mathcal{F},\{\mathcal{F}_t\}_{t\in[0,T]},\mathbb{P})$
be a filtered probability space carrying an $m$-dimensional Brownian
motion
\begin{equation*}
    \mathbf{W}_t=(W_t^1,\ldots,W_t^m)^\top,
\end{equation*}
and let
\begin{equation*}
    \mathcal{H}
    =
    L^2([0,T];\mathbb{R}^m)
\end{equation*}
with inner product
\begin{equation*}
    \langle h_1,h_2\rangle_{\mathcal{H}}
    =
    \int_0^T h_1(s)^\top h_2(s)\,ds.
\end{equation*}
For $h\in\mathcal{H}$, define the Wiener integral
\begin{equation*}
\label{eq:isonormal-process}
    W(h)
    =
    \int_0^T h(s)^\top\,d\mathbf{W}_s
    =
    \sum_{\ell=1}^m
    \int_0^T h^\ell(s)\,dW_s^\ell.
\end{equation*}
The family $\{W(h):h\in\mathcal{H}\}$ is an isonormal Gaussian process
over $\mathcal{H}$; in particular,
\begin{equation*}
    \mathbb{E}[W(h_1)W(h_2)]
    =
    \langle h_1,h_2\rangle_{\mathcal{H}}.
\end{equation*}

Let $\mathcal{S}$ denote the class of smooth cylindrical random variables
of the form
\begin{equation*}
\label{eq:cylindrical-rv}
    F
    =
    f\bigl(W(h_1),\ldots,W(h_n)\bigr),
\end{equation*}
where $h_1,\ldots,h_n\in\mathcal{H}$ and
$f:\mathbb{R}^n\to\mathbb{R}$ is smooth with derivatives of at most
polynomial growth. The Malliavin derivative of $F$ is the
$\mathcal{H}$-valued random variable
\begin{equation*}
\label{eq:malliavin-derivative}
    DF
    =
    \sum_{i=1}^n
    \partial_i f\bigl(W(h_1),\ldots,W(h_n)\bigr)h_i.
\end{equation*}
Equivalently, its value at time $s$ is
\begin{equation*}
\label{eq:malliavin-derivative-time}
    D_sF
    =
    \sum_{i=1}^n
    \partial_i f\bigl(W(h_1),\ldots,W(h_n)\bigr)h_i(s)
    \in\mathbb{R}^m.
\end{equation*}
This definition is independent of the chosen cylindrical
representation. Indeed, after expressing all $h_i$ in an orthonormal
basis of their finite-dimensional span, the assertion follows from the
ordinary chain rule and the nondegeneracy of the resulting Gaussian
vector.

The operator
\begin{equation*}
    D:\mathcal{S}\subset L^2(\Omega)
    \longrightarrow
    L^2(\Omega;\mathcal{H})
\end{equation*}
is closable. Its closed extension has domain
$\mathbb{D}^{1,2}$, defined as the completion of $\mathcal{S}$ under the
norm
\begin{equation*}
\label{eq:malliavin-sobolev-norm}
\begin{aligned}
    \|F\|_{1,2}^2
    &=
    \mathbb{E}[|F|^2]
    +
    \mathbb{E}\bigl[\|DF\|_{\mathcal{H}}^2\bigr]\\
    &=
    \mathbb{E}[|F|^2]
    +
    \mathbb{E}\left[
        \int_0^T |D_sF|^2\,ds
    \right].
\end{aligned}
\end{equation*}
The derivative of a vector- or matrix-valued random variable is defined
componentwise.

The Malliavin derivative satisfies the usual chain rule. If
$F_1,\ldots,F_q\in\mathbb{D}^{1,2}$ and
$\phi\in C^1(\mathbb{R}^q)$ has bounded partial derivatives, then
$\phi(F_1,\ldots,F_q)\in\mathbb{D}^{1,2}$ and
\begin{equation*}
\label{eq:malliavin-chain-rule}
    D\phi(F_1,\ldots,F_q)
    =
    \sum_{i=1}^q
    \partial_i\phi(F_1,\ldots,F_q)\,DF_i.
\end{equation*}
In particular, under the corresponding integrability conditions,
\begin{equation*}
\label{eq:malliavin-product-rule}
    D(FG)=F\,DG+G\,DF.
\end{equation*}

The Cameron--Martin space associated with $\mathbf{W}$ can be identified
with
\begin{equation*}
\begin{split}
    &\mathcal{H}_{\mathrm{CM}}
    =\\
    &\left\{ 
        k_h:[0,T]\to\mathbb{R}^m:
        k_h(t)=\int_0^t h(s)\,ds,\quad h\in\mathcal{H}
    \right\},
\end{split}
\end{equation*}
equipped with the norm
$\|k_h\|_{\mathcal{H}_{\mathrm{CM}}}=\|h\|_{\mathcal{H}}$.
The Malliavin derivative represents differentiation along such
directions. More precisely, for deterministic $h\in\mathcal{H}$, define
\begin{equation*}
\label{eq:directional-malliavin-derivative}
    D_hF
    =
    \langle DF,h\rangle_{\mathcal{H}}
    =
    \int_0^T D_sF\cdot h(s)\,ds.
\end{equation*}
Whenever the corresponding directional derivative exists,
\begin{equation*}
\label{eq:cameron-martin-direction}
    D_hF
    =
    \left.
    \frac{d}{d\varepsilon}
    F(\mathbf{W}+\varepsilon k_h)
    \right|_{\varepsilon=0}.
\end{equation*}

For example, in the one-dimensional case,
\begin{equation*}
    W_t=W(\mathbf{1}_{[0,t]}),
    \quad
    D_sW_t=\mathbf{1}_{[0,t]}(s).
\end{equation*}
More generally, for an $m$-dimensional Brownian motion,
\begin{equation*}
\label{eq:malliavin-derivative-brownian}
    D_s^\ell W_t^k
    =
    \delta_{\ell k}\mathbf{1}_{[0,t]}(s),
    \quad
    1\leq \ell,k\leq m.
\end{equation*}

\paragraph{The Skorokhod integral.}

The divergence operator, also called the Skorokhod integral, is defined
as the adjoint of the Malliavin derivative. A random process
$u\in L^2(\Omega;\mathcal{H})$ belongs to the domain
$\operatorname{Dom}(\delta)$ if there exists a constant $C_u<\infty$
such that
\begin{equation*}
\label{eq:divergence-domain}
    \left|
        \mathbb{E}\bigl[\langle DF,u\rangle_{\mathcal{H}}\bigr]
    \right|
    \leq
    C_u\|F\|_{L^2(\Omega)}
    \quad
    \text{for every }F\in\mathcal{S}.
\end{equation*}
For such a process, $\delta(u)$ is the unique random variable in
$L^2(\Omega)$ satisfying the duality relation
\begin{equation}
\label{eq:malliavin-duality}
    \mathbb{E}\bigl[F\delta(u)\bigr]
    =
    \mathbb{E}\bigl[\langle DF,u\rangle_{\mathcal{H}}\bigr],
    \quad
    F\in\mathbb{D}^{1,2}.
\end{equation}
We use the notation
\begin{equation*}
    \delta(u)
    =
    \int_0^T u_s\cdot\delta\mathbf{W}_s.
\end{equation*}
For deterministic $h\in\mathcal{H}$,
\begin{equation*}
    \delta(h)=W(h).
\end{equation*}
More importantly, if $u$ is adapted and satisfies
\begin{equation*}
    \mathbb{E}\left[\int_0^T |u_s|^2\,ds\right]<\infty,
\end{equation*}
then $u\in\operatorname{Dom}(\delta)$ and its Skorokhod integral
coincides with its It\^o integral:
\begin{equation*}
\label{eq:skorokhod-equals-ito}
    \delta(u)
    =
    \int_0^T u_s\cdot d\mathbf{W}_s.
\end{equation*}
Consequently, for an adapted square-integrable process $u$ and a random
variable $F\in\mathbb{D}^{1,2}$ for which the following expressions are
integrable,
\begin{equation}
\label{eq:malliavin-integration-by-parts}
    \mathbb{E}\left[
        \int_0^T D_sF\cdot u_s\,ds
    \right]
    =
    \mathbb{E}\left[
        F\int_0^T u_s\cdot d\mathbf{W}_s
    \right].
\end{equation}
Equation \eqref{eq:malliavin-integration-by-parts} is the
integration-by-parts identity used below to transfer derivatives of a
functional to stochastic-integral weights.

\paragraph{Malliavin derivative of an SDE flow.}

We finally record the relation between the Malliavin derivative and the
tangent process introduced in Section~\ref{sec:tangent-process}. For
$t\leq s\leq T$, let $\mathbf{X}_s^{t,\mathbf{x}}$ solve
\eqref{eq:sde} with initial value $\mathbf{X}_t = \mathbf{x}$. Under the regularity assumptions imposed in
Section~\ref{sec:fk}, each component of
$\mathbf{X}_s^{t,\mathbf{x}}$ belongs to $\mathbb{D}^{1,2}$.

For $t\leq r\leq s$, the Malliavin derivative
$D_r\mathbf{X}_s^{t,\mathbf{x}}\in\mathbb{R}^{d\times m}$ satisfies
\begin{equation}
\label{eq:malliavin-derivative-sde}
\resizebox{0.88\columnwidth}{!}{$\displaystyle
\begin{split}
    D_r\mathbf{X}_s^{t,\mathbf{x}}
    ={}&
    \boldsymbol{\sigma}
    \bigl(\mathbf{X}_r^{t,\mathbf{x}},r\bigr)
    +
    \int_r^s
    \nabla_{\mathbf{x}}\mathbf{f}
    \bigl(\mathbf{X}_q^{t,\mathbf{x}},q\bigr)
    D_r\mathbf{X}_q^{t,\mathbf{x}}\,dq
    \\
    &+
    \sum_{k=1}^m
    \int_r^s
    \nabla_{\mathbf{x}}\boldsymbol{\sigma}_{\cdot k}
    \bigl(\mathbf{X}_q^{t,\mathbf{x}},q\bigr)
    D_r\mathbf{X}_q^{t,\mathbf{x}}\,dW_q^k.
\end{split}
$}
\end{equation}
Moreover,
\begin{equation*}
    D_r\mathbf{X}_s^{t,\mathbf{x}}=0
    \quad\text{if }r>s.
\end{equation*}
When the Jacobian flow is invertible, the solution of
\eqref{eq:malliavin-derivative-sde} can be written as
\begin{equation*}
\label{eq:malliavin-tangent-relation}
    D_r\mathbf{X}_s^{t,\mathbf{x}}
    =
    \mathbf{J}_{t,s}^{\mathbf{x}}
    \bigl(\mathbf{J}_{t,r}^{\mathbf{x}}\bigr)^{-1}
    \boldsymbol{\sigma}
    \bigl(\mathbf{X}_r^{t,\mathbf{x}},r\bigr),
    \quad
    t\leq r\leq s.
\end{equation*}
In the additive identity-noise setting
$\boldsymbol{\sigma}\equiv\mathbf{I}_d$, this reduces to
\begin{equation*}
\label{eq:malliavin-tangent-identity-noise}
    D_r\mathbf{X}_s^{t,\mathbf{x}}
    =
    \mathbf{J}_{t,s}^{\mathbf{x}}
    \bigl(\mathbf{J}_{t,r}^{\mathbf{x}}\bigr)^{-1},
    \quad
    t\leq r\leq s.
\end{equation*}

\section{Proof of theorems}
\label{app:proofs}

\subsection{Proof of Proposition \ref{prop:reverse-fk-density}}
\label{sec: FP_FK}

Fix $t\in(0,T]$ and define $u(s,\mathbf{y})=p(t-s,\mathbf{y})$ for $s\in[0,t]$. Then $u(0,\mathbf{y})=p(t,\mathbf{y})$ and $u(t,\mathbf{y})=p_0(\mathbf{y})$. Using
 \eqref{eq:fpe-nondiv}, we have that $u$ solves the backward Kolmogorov equation
\begin{equation*}
\begin{aligned}
    &\partial_s u(s,\mathbf{y}) + \boldsymbol{\beta}(t-s,\mathbf{y})\cdot\nabla_{\mathbf{y}}u(s,\mathbf{y}) \\
    &+ \tfrac{1}{2}\,\mathbf{a}(t-s,\mathbf{y}):\nabla^2_{\mathbf{y}}u(s,\mathbf{y}) \\
    &+ c(t-s,\mathbf{y})\,u(s,\mathbf{y})
    = 0,
    \end{aligned}
\end{equation*}
for $s\in[0,t)$, with terminal condition $u(t,\cdot)=p_0$. 
The Feynman--Kac formula \citep{karatzas1991brownian}[Section~5.7] therefore gives
\begin{equation*}
\begin{aligned}
    u(0,\mathbf{x}) &= \mathbb{E}\left[p_0\bigl(\mathbf{Y}_t^{t,\mathbf{x}}\bigr)\exp\left(\int_0^t c\bigl(t-s,\mathbf{Y}_s^{t,\mathbf{x}}\bigr)\,ds\right)\right] \\
    &= \mathbb{E}\bigl[G_t^{t,\mathbf{x}}\bigr].
    \end{aligned}
\end{equation*}
Since $u(0,\mathbf{x})=p(t,\mathbf{x})$ and using \eqref{eq:fk-payoff-general}, we have this is \eqref{eq:rep-fpe}.

\subsection{Proof of Theorem \ref{thm:BEL}}
\label{sec:proof_thm_BEL}
\begin{proof}
Fix $t\in(0,T]$, $\mathbf{x}\in\mathbb{R}^d$, and a direction
$\mathbf{v}\in\mathbb{R}^d$. To simplify notation, throughout the
proof we write
\begin{align*}
&\mathbf{Y}_s=\mathbf{Y}^{t,\mathbf{x}}_s,
\quad
&&\mathbf{J}_s=\mathbf{J}^{t,\mathbf{x}}_s,
\quad\\
&\boldsymbol{\sigma}^t_s
    =\boldsymbol{\sigma}(t-s,\mathbf{Y}_s),
\quad
&&\mathbf{a}^t_s
    =\mathbf{a}(t-s,\mathbf{Y}_s),
\end{align*}
and define
\begin{equation*}
C_t=\int_0^t c(t-s,\mathbf{Y}_s)\,ds,
\qquad
G=p_0(\mathbf{Y}_t)e^{C_t}.
\end{equation*}
We first assume that $p_0\in C_b^1(\mathbb{R}^d)$. The general case
follows by approximation at the end of the proof.

Differentiating the Feynman--Kac representation with respect to the
initial point in the direction $\mathbf{v}$ gives
\begin{equation}\label{eq:proof-bel-pathwise}
\resizebox{0.881\columnwidth}{!}{$\displaystyle
\begin{split}
\nabla_{\mathbf{x}}p(t,\mathbf{x})\cdot\mathbf{v}
={}&
\mathbb{E}\!\left[
e^{C_t}\,
\nabla p_0(\mathbf{Y}_t)^\top\mathbf{J}_t\mathbf{v}
\right] \\
&+
\mathbb{E}\!\left[
G\int_0^t
\nabla_{\mathbf{x}}c(t-s,\mathbf{Y}_s)^\top
\mathbf{J}_s\mathbf{v}\,ds
\right].
\end{split}
$}
\end{equation}
The second term already has the desired form. We use Malliavin
integration by parts to remove $\nabla p_0$ from the first term.

Let
\begin{equation*}
(\boldsymbol{\sigma}^t_s)^+
=
(\boldsymbol{\sigma}^t_s)^\top(\mathbf{a}^t_s)^{-1},
\end{equation*}
so that
\begin{equation*}
\boldsymbol{\sigma}^t_s(\boldsymbol{\sigma}^t_s)^+
=\mathbf{I}_d.
\end{equation*}
Define the adapted $\mathbb{R}^m$-valued process
\begin{equation}\label{eq:proof-bel-u}
\mathbf{u}^{\mathbf{v}}_r
=
\frac{1}{t}
(\boldsymbol{\sigma}^t_r)^+
\mathbf{J}_r\mathbf{v}\,
\mathbf{1}_{[0,t]}(r).
\end{equation}
Under the assumptions of the theorem,
$\mathbf{u}^{\mathbf{v}}\in\operatorname{Dom}(\delta)$, and since it
is adapted,
\begin{equation}\label{eq:proof-bel-delta}
\begin{split}
\delta(\mathbf{u}^{\mathbf{v}})
&=
\frac{1}{t}\int_0^t
\bigl((\boldsymbol{\sigma}^t_r)^+
      \mathbf{J}_r\mathbf{v}\bigr)^\top
\,d\mathbf{W}_r \\
&=
\mathbf{v}^\top
\left[
\frac{1}{t}\int_0^t
\mathbf{J}_r^\top(\mathbf{a}^t_r)^{-1}
\boldsymbol{\sigma}^t_r\,d\mathbf{W}_r
\right].
\end{split}
\end{equation}

For $0\leq r\leq s\leq t$, the relation between the Malliavin
derivative and the tangent process gives
\begin{equation*}
D_r\mathbf{Y}_s
=
\mathbf{J}_s\mathbf{J}_r^{-1}
\boldsymbol{\sigma}^t_r.
\end{equation*}
Consequently,
\begin{align}
\left\langle
D_r p_0(\mathbf{Y}_t),
\mathbf{u}^{\mathbf{v}}_r
\right\rangle_{\mathbb{R}^m}
&=
\frac{1}{t}
\nabla p_0(\mathbf{Y}_t)^\top
D_r\mathbf{Y}_t
(\boldsymbol{\sigma}^t_r)^+
\mathbf{J}_r\mathbf{v}
\nonumber\\
&=
\frac{1}{t}
\nabla p_0(\mathbf{Y}_t)^\top
\mathbf{J}_t\mathbf{v}.
\label{eq:proof-bel-collapse}
\end{align}

Notice that the right hand of \eqref{eq:proof-bel-collapse} is independent of $r$. Integrating \eqref{eq:proof-bel-collapse} over $r\in[0,t]$ and multiplying $e^{C_t}$ yields
\begin{equation}
    \label{eq:malliavin-intermediate}
    e^{C_t}\nabla p_0(\mathbf{Y}_t)^\top\mathbf{J}_t\mathbf{v}
=
e^{C_t}
\left\langle
D p_0(\mathbf{Y}_t),\mathbf{u}^{\mathbf{v}}
\right\rangle_{\mathcal H}.
\end{equation}

By the Malliavin product rule,
\begin{equation*}
D G
=
e^{C_t}D p_0(\mathbf{Y}_t)+G\,D C_t.
\end{equation*}
Therefore, integrating over $r \in [0, t]$ again, we obtain
\begin{equation*}\label{eq:proof-bel-product}
e^{C_t}
\left\langle
D p_0(\mathbf{Y}_t),\mathbf{u}^{\mathbf{v}}
\right\rangle_{\mathcal H}
=
\left\langle DG,\mathbf{u}^{\mathbf{v}}\right\rangle_{\mathcal H}
-
G\left\langle DC_t,\mathbf{u}^{\mathbf{v}}\right\rangle_{\mathcal H}.
\end{equation*}
Take expectations, apply equation \eqref{eq:malliavin-intermediate}, and apply the duality relation between $D$
and $\delta$ in equation \eqref{eq:malliavin-duality}. Now we have
\begin{equation}\label{eq:proof-bel-duality}
\resizebox{0.88\columnwidth}{!}{$\displaystyle
\mathbb{E}\!\left[
e^{C_t}\nabla p_0(\mathbf{Y}_t)^\top
\mathbf{J}_t\mathbf{v}
\right]
=
\mathbb{E}\!\left[G\delta(\mathbf{u}^{\mathbf{v}})\right]
-
\mathbb{E}\!\left[
G\left\langle DC_t,\mathbf{u}^{\mathbf{v}}\right\rangle_{\mathcal H}
\right].
$}
\end{equation}

The first expectation on the right hand side is given by equation \eqref{eq:proof-bel-delta}. It remains to calculate the second expectation. For $r\leq t$,
the Malliavin chain rule gives
\begin{equation*}
D_rC_t
=
\int_r^t
(D_r\mathbf{Y}_s)^\top
\nabla_{\mathbf{x}}c(t-s,\mathbf{Y}_s)\,ds.
\end{equation*}
Using the expression for $D_r\mathbf{Y}_s$ and
$\boldsymbol{\sigma}^t_r(\boldsymbol{\sigma}^t_r)^+
=\mathbf{I}_d$, we obtain
\begin{equation*}
\begin{aligned}
&\left\langle D_rC_t,\mathbf{u}^{\mathbf{v}}_r
\right\rangle_{\mathbb{R}^m}
\\
&\quad=
\frac{1}{t}\int_r^t
\nabla_{\mathbf{x}}c(t-s,\mathbf{Y}_s)^\top
\mathbf{J}_s\mathbf{J}_r^{-1}
\boldsymbol{\sigma}^t_r
(\boldsymbol{\sigma}^t_r)^+
\mathbf{J}_r\mathbf{v}\,ds
\\
&\quad=
\frac{1}{t}\int_r^t
\nabla_{\mathbf{x}}c(t-s,\mathbf{Y}_s)^\top
\mathbf{J}_s\mathbf{v}\,ds.
\end{aligned}
\end{equation*}
Fubini's theorem then yields
\begin{equation}\label{eq:proof-bel-potential}
\small
\begin{split}
\left\langle DC_t,\mathbf{u}^{\mathbf{v}}\right\rangle_{\mathcal H}
&=
\frac{1}{t}\int_0^t\int_r^t
\nabla_{\mathbf{x}}c(t-s,\mathbf{Y}_s)^\top
\mathbf{J}_s\mathbf{v}\,ds\,dr \\
&= \frac{1}{t}\int_0^t\int_0^s
\nabla_{\mathbf{x}}c(t-s,\mathbf{Y}_s)^\top
\mathbf{J}_s\mathbf{v}\,dr\,ds \\
&=
\int_0^t\frac{s}{t}
\nabla_{\mathbf{x}}c(t-s,\mathbf{Y}_s)^\top
\mathbf{J}_s\mathbf{v}\,ds.
\end{split}
\end{equation}

Substituting \eqref{eq:proof-bel-duality} and
\eqref{eq:proof-bel-potential} into
\eqref{eq:proof-bel-pathwise}, the contribution
$s/t$ is subtracted from the pathwise derivative of the potential.
Hence
\begin{align*}
&\nabla_{\mathbf{x}}p(t,\mathbf{x})\cdot\mathbf{v}\\
&=
\mathbb{E}\!\Bigg[
G\Bigg\{
\delta(\mathbf{u}^{\mathbf{v}})+
\int_0^t\frac{t-s}{t}
\nabla_{\mathbf{x}}c(t-s,\mathbf{Y}_s)^\top
\mathbf{J}_s\mathbf{v}\,ds
\Bigg\}
\Bigg].
\end{align*}
Using \eqref{eq:proof-bel-delta}, the expression inside braces equals
$\mathbf{v}^\top\mathbf{H}^{t,\mathbf{x}}_t$. We have therefore shown
that
\begin{equation*}
\nabla_{\mathbf{x}}p(t,\mathbf{x})\cdot\mathbf{v}
=
\mathbf{v}^\top
\mathbb{E}\!\left[
G^{t,\mathbf{x}}_t\mathbf{H}^{t,\mathbf{x}}_t
\right].
\end{equation*}
Since $\mathbf{v}\in\mathbb{R}^d$ is arbitrary,
\begin{equation*}
\nabla_{\mathbf{x}}p(t,\mathbf{x})
=
\mathbb{E}\!\left[
G^{t,\mathbf{x}}_t\mathbf{H}^{t,\mathbf{x}}_t
\right].
\end{equation*}

Finally, for a bounded continuous $p_0$, let
$p_0^{(n)}$ be a bounded smooth mollification of $p_0$. The identity
above applies to each $p_0^{(n)}$. Uniform ellipticity and the moment
bounds for $\mathbf{J}$ and $(\boldsymbol{\sigma}^t_s)^+$, together
with the assumed integrability of $e^{C_t}$ and
$\mathbf{H}^{t,\mathbf{x}}_t$, allow passage to the limit by the
standard localization and dominated-convergence argument for the
Feynman--Kac semigroup. Thus the formula remains valid without
assuming that $p_0$ is differentiable.
  
\end{proof}

\subsection{Proof of Theorem \ref{prop:score-asymptotic-variance}}
\label{sec:proof_variance}
For the sake of simplicity, write $p=\mathbb{E}[G]$, $\mathbf{g}=\mathbb{E}[G\mathbf{H}]$, and $\mathbf{s}=\mathbf{g}/p$. By the strong law of large numbers, we have almost surely convergence
\begin{equation*}
    \widehat p_N\longrightarrow p,
    \quad
    \widehat{\mathbf{g}}_N\longrightarrow\mathbf{g}
    \,.
\end{equation*}
Since $p>0$, the continuous mapping theorem gives the almost surely convergence
\begin{equation*}
    \widehat{\mathbf{s}}_N
    =
    \frac{\widehat{\mathbf{g}}_N}{\widehat p_N}
    \longrightarrow
    \frac{\mathbf{g}}{p}
    =
    \mathbf{s} \,.
\end{equation*}

To obtain the asymptotic distribution, define the
$(d+1)$-dimensional random vector
\begin{equation*}
    \mathbf{Z}
    =
    \begin{pmatrix}
        G\\
        G\mathbf{H}
    \end{pmatrix},
    \quad
    \boldsymbol{\mu}
    =
    \mathbb{E}[\mathbf{Z}]
    =
    \begin{pmatrix}
        p\\
        \mathbf{g}
    \end{pmatrix}.
\end{equation*}
The moment assumptions imply
that $\mathbf{Z}$ has a finite covariance matrix, denoted by
$\boldsymbol{\Omega}$. The multivariate central limit theorem gives
\begin{equation*}
\label{eq:joint-density-gradient-clt}
    \sqrt{N}
    \left[
        \begin{pmatrix}
            \widehat p_N\\
            \widehat{\mathbf{g}}_N
        \end{pmatrix}
        -
        \begin{pmatrix}
            p\\
            \mathbf{g}
        \end{pmatrix}
    \right]
    \xrightarrow{\mathrm{d}}
    \mathcal{N}
    \left(
        \mathbf{0},
        \boldsymbol{\Omega}
    \right).
\end{equation*}

It remains to use the Delta method to estimate $\bm{\Omega}$. Consider the smooth map $\Phi:
    (0,\infty)\times\mathbb{R}^d
    \longrightarrow
    \mathbb{R}^d$ for $\Phi(a,\mathbf{b}) = \mathbf{b}/a$. Its Jacobian matrix at $(p,\mathbf{g})$ is
\begin{equation*}
\label{eq:ratio-map-jacobian}
    D\Phi(p,\mathbf{g})
    =
    \begin{pmatrix}
        -\mathbf{g}/p^2 & \mathbf{I}_d/p
    \end{pmatrix}
    =
    \frac{1}{p}
    \begin{pmatrix}
        -\mathbf{s} & \mathbf{I}_d
    \end{pmatrix}.
\end{equation*}
The multivariate delta method therefore yields
\begin{equation*}
    \sqrt{N}
    \left(
        \widehat{\mathbf{s}}_N-\mathbf{s}
    \right)
    \xrightarrow{\mathrm{d}}
    \mathcal{N}
    \left(
        \mathbf{0},
        D\Phi(p,\mathbf{g})
        \boldsymbol{\Omega}
        D\Phi(p,\mathbf{g})^\top
    \right).
\end{equation*}

The covariance matrix can be further simplified. Notice that
\begin{align*}
    D\Phi(p,\mathbf{g})
    \left(
        \mathbf{Z}-\boldsymbol{\mu}
    \right)
    &=
    \frac{1}{p}
    \left[
        G\mathbf{H}-\mathbf{g}
        -
        \mathbf{s}(G-p)
    \right]
    \\
    &=
    \frac{G}{p}
    \left(
        \mathbf{H}-\mathbf{s}
    \right),
\end{align*}
where we used $\mathbf{g}=p\mathbf{s}$. Its expectation is zero because
\begin{equation*}
    \mathbb{E}\left[
        G(\mathbf{H}-\mathbf{s})
    \right]
    =
    \mathbf{g}-p\mathbf{s}
    =
    \mathbf{0}.
\end{equation*}
Consequently,
\begin{equation*}
    D\Phi(p,\mathbf{g})
    \boldsymbol{\Omega}
    D\Phi(p,\mathbf{g})^\top
    =
    \frac{
        \mathbb{E}\left[
            G^2
            (\mathbf{H}-\mathbf{s})
            (\mathbf{H}-\mathbf{s})^\top
        \right]
    }{
        p^2
    },
\end{equation*}
which is the expression of $\boldsymbol{\Sigma}(t,\mathbf{x}) $ in equation \eqref{eq:score-clt}. 

\subsection{Proof of Proposition \ref{prop:ou-overlap-growth} and discussions} 
\label{sec:ou-variance}
\label{sec:OU_proof}

\begin{proof}
The solution of \eqref{eq:reverse-ou} is
\begin{equation*}
    \mathbf{Y}_t
    =
    e^{\lambda t}\mathbf{x}
    +
    \sigma
    \int_0^t
    e^{\lambda(t-s)}\,d\mathbf{W}_s.
\end{equation*}
Therefore,
\begin{equation*}
\label{eq:reverse-ou-law}
    \mathbf{Y}_t
    \sim
    \mathcal{N}
    \left(
        e^{\lambda t}\mathbf{x},
        \nu_t\mathbf{I}_d
    \right).
\end{equation*}
For a Gaussian random variable
$\mathbf{Z}\sim\mathcal{N}(\mathbf{m},\nu\mathbf{I}_d)$ and any
$\alpha>0$,
\begin{equation}
\label{eq:gaussian-quadratic-integral}
\resizebox{0.88\columnwidth}{!}{$\displaystyle
    \mathbb{E}
    \left[
        \exp\left(
            -\frac{\alpha}{2}\|\mathbf{Z}\|^2
        \right)
    \right]
    =
    (1+\alpha\nu)^{-d/2}
    \exp\left(
        -\frac{\alpha\|\mathbf{m}\|^2}
              {2(1+\alpha\nu)}
    \right).
$}
\end{equation}
Applying \eqref{eq:gaussian-quadratic-integral} with
$\alpha=\tau_0^{-2}$ gives
\begin{equation}
\label{eq:ou-first-payoff-moment}
\begin{split}
    \mathbb{E}[p_0(\mathbf{Y}_t)]
    ={}&
    \frac{1}{(2\pi\tau_0^2)^{d/2}}
    (1+r_t)^{-d/2}
    \\
    &\times
    \exp\left(
        -\frac{
            e^{2\lambda t}\|\mathbf{x}\|^2
        }{
            2(\tau_0^2+\nu_t)
        }
    \right).
\end{split}
\end{equation}
Similarly, applying it with $\alpha=2\tau_0^{-2}$ gives
\begin{equation}
\label{eq:ou-second-payoff-moment}
\begin{split}
    \mathbb{E}[p_0(\mathbf{Y}_t)^2]
    ={}&
    \frac{1}{(2\pi\tau_0^2)^d}
    (1+2r_t)^{-d/2}
    \\
    &\times
    \exp\left(
        -\frac{
            e^{2\lambda t}\|\mathbf{x}\|^2
        }{
            \tau_0^2+2\nu_t
        }
    \right).
\end{split}
\end{equation}
Notice that we have $G = p_0(\mathbf{Y}_t)$ for OU process. Substituting \eqref{eq:ou-first-payoff-moment} and
\eqref{eq:ou-second-payoff-moment} into
\eqref{eq:overlap-factor} yields
\eqref{eq:ou-R2-exact}.
\end{proof}

Proposition \ref{prop:ou-overlap-growth} implies that $R_2(t, \mathbf{x})$ is exponential in both time and dimension. Notice that 
\begin{equation*}
    \frac{1}{\tau_0^2+\nu_t}
    -
    \frac{1}{\tau_0^2+2\nu_t}
    >0.
\end{equation*}
So the exponential term in \eqref{eq:ou-R2-exact} is at least one. Moreover, for $r_t>0$,
\begin{equation*}
    (1+r_t)^2-(1+2r_t)=r_t^2>0,
\end{equation*}
which implies $\kappa_t > 1$ for $\kappa_t = (1 + r_t)/ \sqrt{1 + 2 r_t}$, which gives $R_2(t,x) > \kappa_t^d$. As $t \rightarrow \infty$, we have  
\begin{equation*}
    r_t
    =
    \frac{\sigma^2}{2\lambda\tau_0^2}
    e^{2\lambda t}
    \bigl(1+o(1)\bigr).
\end{equation*}
Hence
\begin{equation*}
    \kappa_t = \frac{1+r_t}{\sqrt{1+2r_t}}
    =
    \frac{\sigma}
         {2\sqrt{\lambda}\tau_0}
    e^{\lambda t}
    \bigl(1+o(1)\bigr),
\end{equation*}
and $
R_2(t, \mathbf{x}) \approx O(1) e^{\lambda d t} $. Hence the overlap factor grows exponentially in time, with an asymptotic exponential rate proportional to the dimension.

Proposition~\ref{prop:ou-overlap-growth} shows that 
the variance amplification is mainly caused by the increasingly poor overlap
between the unstable reverse process and the initial density. In high dimensions, this overlap deterioration means the BEL score estimate $\hat{\mathbf{s}}$ is only useful for very short $t$.

\section{Experiments}\label{app:experiments}

\subsection{Mixing rate and expansion rate of ABC-like flows}\label{app:mixing-vs-expansion}
This experiment considers the ABC-like flow 
\begin{equation*}
\mathrm{d} x_i = [\sin(X_{i-1}) + \cos(X_{i+1})] dt + \sigma dW^i_t
\end{equation*}
as described in equation \eqref{eq:ABC}. We compare the exponential rate of mixing and the exponential rate of variance expansion. For the mixing rate, we examine the speed of convergence of the expectation of the observable
\begin{equation*}
\phi(\mathbf{x}) = \frac{1}{d} \sum_{i = 1}^d \cos(\mathbf{x}_i) + \frac{1}{2} \cos(\sum_{i = 1}^d \mathbf{x}_i) \,.
\end{equation*}
The expectation of $\phi$ with respect to the uniform distribution $u$ on the torus is zero. For the rate of expansion, we use Monte Carlo simulation to estimate the matrix norm (2-norm) of $\mathbf{H}$ in equation \eqref{eq:div-free-malliavin-weight}. It\^o isometry gives 
\begin{equation*}
\mathbb{E}[\|\mathbf{H}_t \|^2]  =\frac{1}{\sigma^2 t^2}\mathbb{E} \int_0^t ||\mathbf{J}_s||_F^2 ds.
\end{equation*}

\begin{figure}
    \centering
    \includegraphics[width=\linewidth]{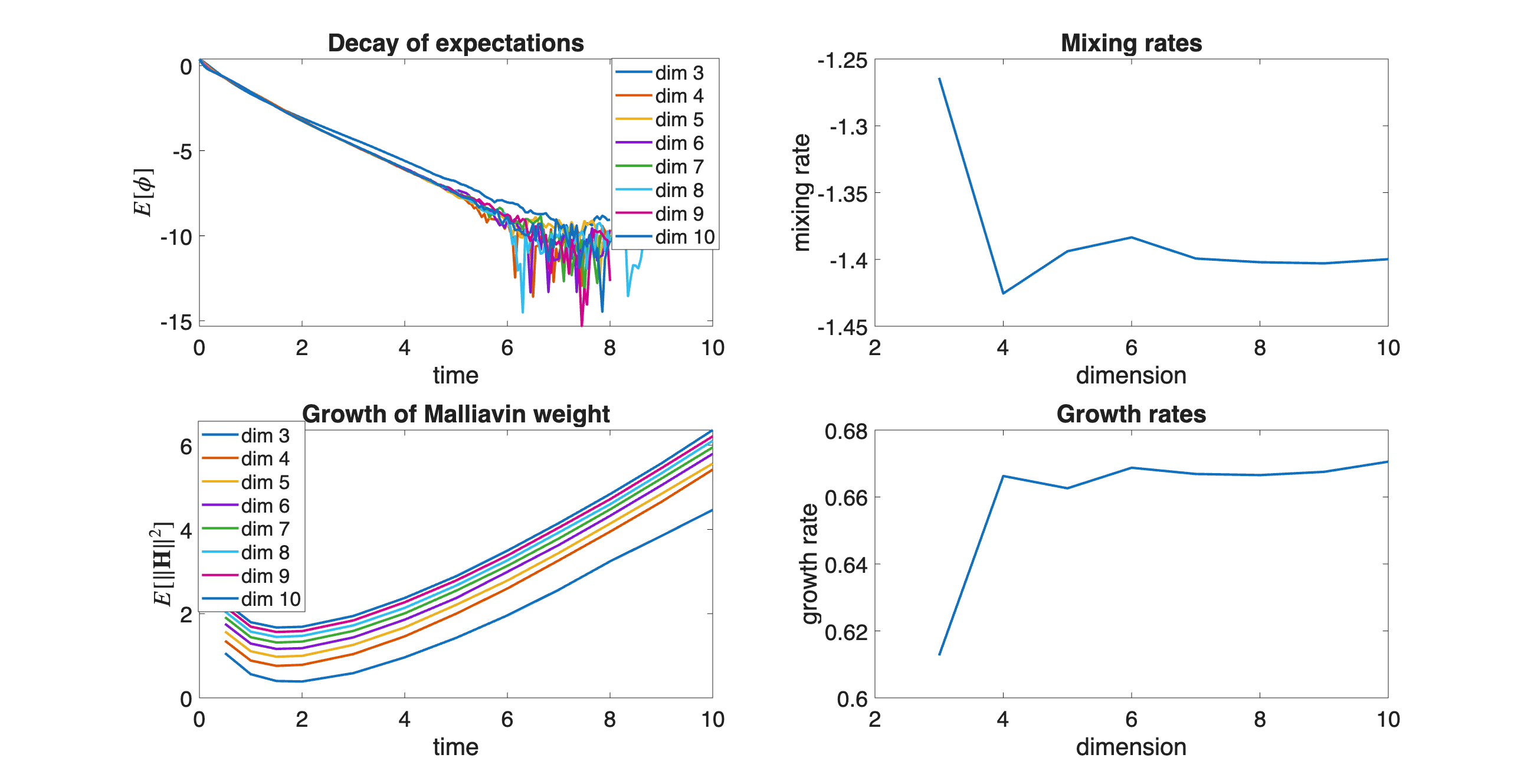}
    \caption{Top: Decay of $\mathbb{E}[\phi(\mathbf{X}_t)]$ in log-linear plot and estimated mixing rate. Bottom: Expansion of $\mathbb{E}[\|\mathbf{H}_t \|^2]$ and estimated expansion rate.}
    \label{fig:mixing}
\end{figure}

In this numerical experiment, we use Monte Carlo simulation to compute $\mathbb{E}[\phi(\mathbf{X}_t)]$ and $\mathbb{E}[\|\mathbf{H}_t \|^2]$ for $d = 3, 4, \cdots, 10$. The Monte Carlo simulation of $\mathbb{E}[\phi(\mathbf{X}_t)]$ uses $10^9$ independent trajectories. The Monte Carlo simulation of $\mathbb{E}[\|\mathbf{H}_t \|^2]$ uses $65536$ trajectories. The results are shown in Figure \ref{fig:mixing}. As seen in Figure \ref{fig:mixing}, both the mixing rate and the expansion rate are stable with increasing dimension. In addition, possibly because of the additive noise, the mixing rate is greater than the expansion rate of $\mathbb{E}[\|\mathbf{H}_t \|^2]$. Let $\varepsilon$ be the targeted approximation error of $q_t^{\mathbf{x}}(\mathbf{y})$ with respect to the invariant probability measure $u$. Let $r_{m}$ be the mixing rate. Then the time horizon of the simulation is $t_m = |\log(\varepsilon)| r_m^{-1}$. Let $r_e$ be the expansion rate of $\mathbb{E}[\|\mathbf{H}_t \|^2]$. Then at $t_m$, we have
\[
\mathbb{E}[\|\mathbf{H}_t \|^2] \approx \varepsilon^{-r_e/r_m} \,.
\]
Therefore, higher ratio of the mixing rate to the expansion rate significantly reduces the number of samples needed. 

\subsection{Ablation Study: number of collocation points}
\label{app:bel_number_spatial}

{\bf Model selection and training implementation. }
The initial distribution is a four-component wrapped Gaussian mixture model (GMM4) on $\mathbb T^3=[-\pi,\pi)^3$, defined by
\begin{equation*}
p_0(x) = \sum_{k=1}^{4}w_k\sum_{n\in\mathbb Z^3}\mathcal N(x+2\pi n;\mu_k,s_k^2I_3),
\end{equation*}
where
\begin{equation*}
(w_1,w_2,w_3,w_4)=(0.35,0.25,0.25,0.15),
\end{equation*}
\begin{equation*}
\begin{aligned}
\mu_1&=(-1.60,-1.20,0.50),&
\mu_2&=(1.40,-1.90,-1.10),\\
\mu_3&=(0.30,1.70,1.80),&
\mu_4&=(-2.40,2.20,-2.30),
\end{aligned}
\end{equation*}
and
\begin{equation*}
(s_1,s_2,s_3,s_4)=(0.45,0.60,0.35,0.70).
\end{equation*}



\begin{figure}[h!]
    \centering
    \includegraphics[width=0.9\linewidth]{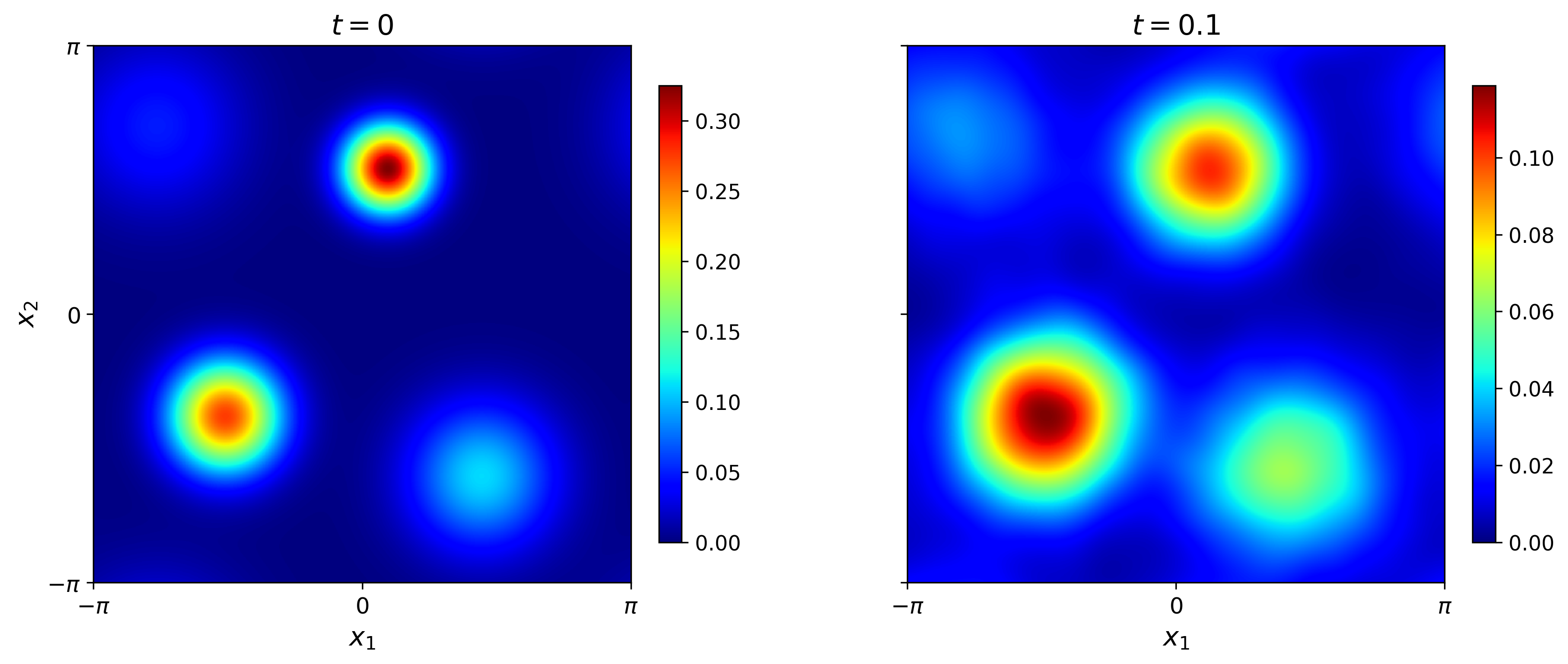}
    \caption{$(x_1,x_2)$-marginal densities of the ABC-flow distribution at $t=0$ (left) and $t=0.1$ (right). Each panel uses an independent color scale.}
    \label{fig:marginal_density_GMM4}
\end{figure}

Figure~\ref{fig:marginal_density_GMM4} shows the $(x_1,x_2)$-marginal densities
\begin{equation*}
p_t^{(1,2)}(x_1,x_2)=\int_{-\pi}^{\pi}p_t(x_1,x_2,x_3)\,dx_3
\end{equation*}
at $t=0$ and $t=0.1$. The four projected modes remain visible at $t=0.1$, with broader and lower peaks. These marginals provide qualitative context for the score-error distributions on the slice $x_3=0$ examined below.

Both networks use five hidden layers of width 128, swish activations, and three outputs, with periodic inputs $(\sin x, \cos x, t)$. All ablation experiments use 40,000 Adam steps with batch size 1024 per active loss and cosine learning-rate decay from $10^{-3}$ to $5 \times 10^{-5}$.

In the score-PINN method, the PDE residual is evaluated at 200,000 space-time collocation points, constructed from 20,000 uniformly distributed spatial points on $\mathbb T^3$, each paired with 10 time coordinates selected from $[0,3]$ according to the nonuniform temporal sampling strategy described in Section~\ref{sec:practical-implementation}.

For the ISM method, we simulate 20,000 trajectories of the forward stochastic dynamics, with initial conditions sampled from $p_0$. Each trajectory is recorded at 300 equally spaced time points $t_k=0.01k$, $k=1,\ldots,300$. During training, each ISM sample is obtained by randomly selecting a trajectory and a time index, and taking the corresponding state-time pair $(X_{t_k}^{(i)},t_k)$, $i=1,\ldots,20000$.

For the BEL anchor, labels are computed from Monte Carlo estimates of the density and its spatial gradient using the ratio in \eqref{eq:mc-score-estimator}. We retain only labels satisfying
\begin{equation*}
    \frac{\widehat p(x,t)}
{\operatorname{SE}[\widehat p(x,t)]}>10,
\end{equation*}
where \(\operatorname{SE}\) denotes the Monte Carlo standard error. The BEL supervised loss uses weights proportional to the estimated density, normalized to have unit mean over the retained training observations. The validation root-mean-square error (RMSE) is the square root of the density-weighted mean squared discrepancy between the predicted scores and the retained BEL labels.

For positive $N_{\mathrm{BEL}}, 15 \%$ of spatial points are reserved for validation (rounded, with at least one held-out point), and the remainder is used for training. All times for each point remain together. The split precedes label filtering. The checkpoint minimizing the density-weighted BEL validation RMSE is selected from evaluations every 500 steps and at the final step; without anchors, the final checkpoint is used. The ablation of collocation-point number uses nested spatial subsets, with the training-validation split recomputed at each budget.

{\bf Effect of the number of collocation points. }

Table~\ref{tab:bel_number_error} reports the relative $L^2$ score errors corresponding to Figure~\ref{fig:bel_number_error} in the main text, providing the numerical values for different numbers of BEL collocation points. The errors for ISM and score-PINN are evaluated using their respective numerical references. 

\begin{table}[h!]
\centering
\begin{tabular}{ccc}
\toprule
$N_{\mathrm{BEL}}$
& ISM + BEL
& Score-PINN + BEL \\
\midrule
20000 & 0.054583 & 0.021677 \\
10000 & 0.054958 & 0.023267 \\
5000  & 0.069788 & 0.023671 \\
2500  & 0.087092 & 0.025239 \\
1250  & 0.116045 & 0.030359 \\
625   & 0.142332 & 0.045409 \\
300   & 0.232723 & 0.055052 \\
150   & 0.255866 & 0.111871 \\
50    & 0.272390 & 0.328210 \\
20    & 0.314583 & 0.617055 \\
10    & 0.321308 & 0.867815 \\
2     & 0.302494 & 0.989438 \\
0     & 0.513801 & 0.999999 \\
\bottomrule
\end{tabular}
\caption{
Relative $L^2$ score errors for ISM and score-PINN
with different numbers of BEL collocation points.
}
\label{tab:bel_number_error}
\end{table}

\begin{figure}[h!]
    \centering
    \includegraphics[width=\linewidth]{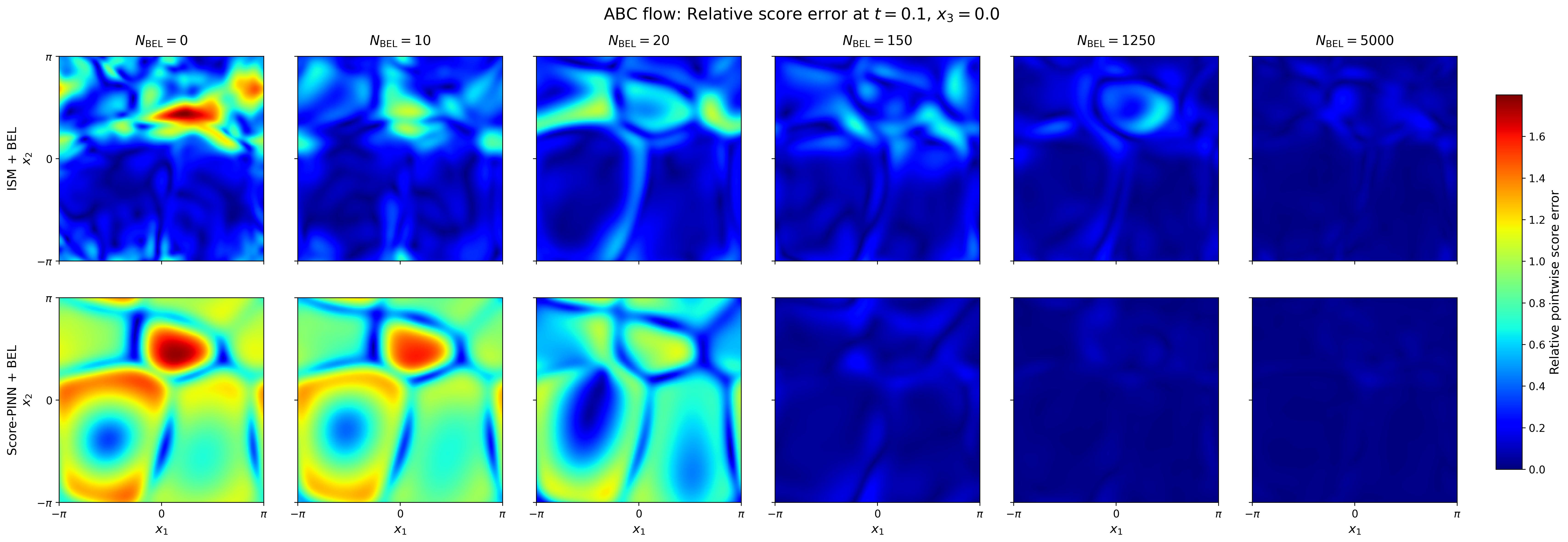}
    \caption{ Spatial distributions of the relative pointwise score errors for ISM with BEL anchor (top row) and score-PINN with BEL supervision (bottom row) with $N_{\mathrm{BEL}}=0,10,20,150,1250,5000$. The errors are evaluated at $t=0.1$ on the spatial slice $x_3=0$, using the related method-specific numerical reference.}
    \label{fig:bel_heatmap_L2_error}
\end{figure}

To illustrate the effect of reducing BEL supervision, we examine the spatial distributions of the relative pointwise score errors for ISM with BEL supervision and score-PINN with BEL supervision at $t=0.1$, where the probability density retains pronounced nonuniform features inherited from the initial distribution. Figure~\ref{fig:bel_heatmap_L2_error} presents the errors on the two-dimensional slice $(x_1,x_2)$ at $x_3=0$ of $\mathbb T^3$, using the method-specific numerical references $s_{\mathrm{ref}}^{\mathrm{ISM}}$ and $s_{\mathrm{ref}}^{\mathrm{PINN}}$ and the same color scale across all panels.

As shown in Figure~\ref{fig:bel_heatmap_L2_error}, the two methods exhibit different spatial error patterns. When BEL supervision is absent or very limited, score-PINN exhibits relatively large errors over broad spatial regions, whereas the errors of ISM are more localized in low coverage areas. Comparison with the $t=0.1$ marginal in Figure~\ref{fig:marginal_density_GMM4} shows that pronounced errors appear near the upper projected mode, while the lower-left region generally exhibits smaller errors. As the number of BEL collocation points increases, the spatial errors decrease for both methods, with a more pronounced improvement observed for score-PINN. These observations are consistent with the relative $L^2$ error trends reported in Section~\ref{sec:ablation_study}.

\subsection{Ablation Study: spatial distribution of collocation points} \label{app:bel_strategies_spatial}


In addition to the number of BEL collocation points,
their spatial distribution may also influence the accuracy
of score estimation. To investigate this effect, we compare
five spatial sampling strategies for BEL supervision, namely Uniform, $p_0$, Half, Gradient, and Complement.
Their corresponding sampling densities on $\mathbb T^3$
are defined by
\begin{equation*}
\begin{aligned}
q_{\text{uni}}(x)&=\frac{1}{(2\pi)^3},\\
q_{p_0}(x)&=p_0(x),\\
q_{\text{half}}(x)&=\frac12q_{\text{uni}}(x)+\frac12p_0(x),\\
q_{\text{grad}}(x)&\propto\|\nabla p_0(x)\|_2,\\
q_{\text{comp}}(x)&\propto p_{0,\max}-p_0(x),
\end{aligned}
\end{equation*}
where $p_{0,\max}=\max_{x\in\mathbb T^3}p_0(x)$.
Gradient sampling favors regions where the initial density varies rapidly, while Complement sampling favors regions of lower initial density.

We fix the number of BEL spatial collocation points at $N_{\mathrm{BEL}}=1250$, motivated by the noticeable increase in relative $L^2$ error as $N_{\mathrm{BEL}}$ decreases from 1250 to 625 in Figure~\ref{fig:bel_number_error}. For each strategy, spatial collocation points are sampled according to the corresponding distribution, while the temporal coordinates follow the same sampling procedure, with identical time coordinates assigned to the same point ID across all five strategies.
All other settings, including the score-PINN collocation dataset, ISM trajectory dataset, network architecture, and training parameters, remain unchanged.

The relative $L^2$ score error is evaluated using 50,000 uniformly sampled spatial points on $\mathbb{T}^3$ and is defined by
\begin{equation*}
E_{L^2}^{m}(t)=\sqrt{\frac{\sum_j\|\mathbf{s}_{\boldsymbol{\theta}}^{m}(t,\boldsymbol{x}_j)-\mathbf{s}_{\mathrm{ref}}^{m}(t,\boldsymbol{x}_j)\|_2^2}{\sum_j\|\mathbf{s}_{\mathrm{ref}}^{m}(t,\boldsymbol{x}_j)\|_2^2}}\;,
\end{equation*}
where $m\in\{\mathrm{ISM},\mathrm{PINN}\}$ and $\mathbf{s}_{\mathrm{ref}}^{m}$ denotes the corresponding numerical reference for each method.

\begin{figure}[h!]
    \centering
    \includegraphics[width=\linewidth]{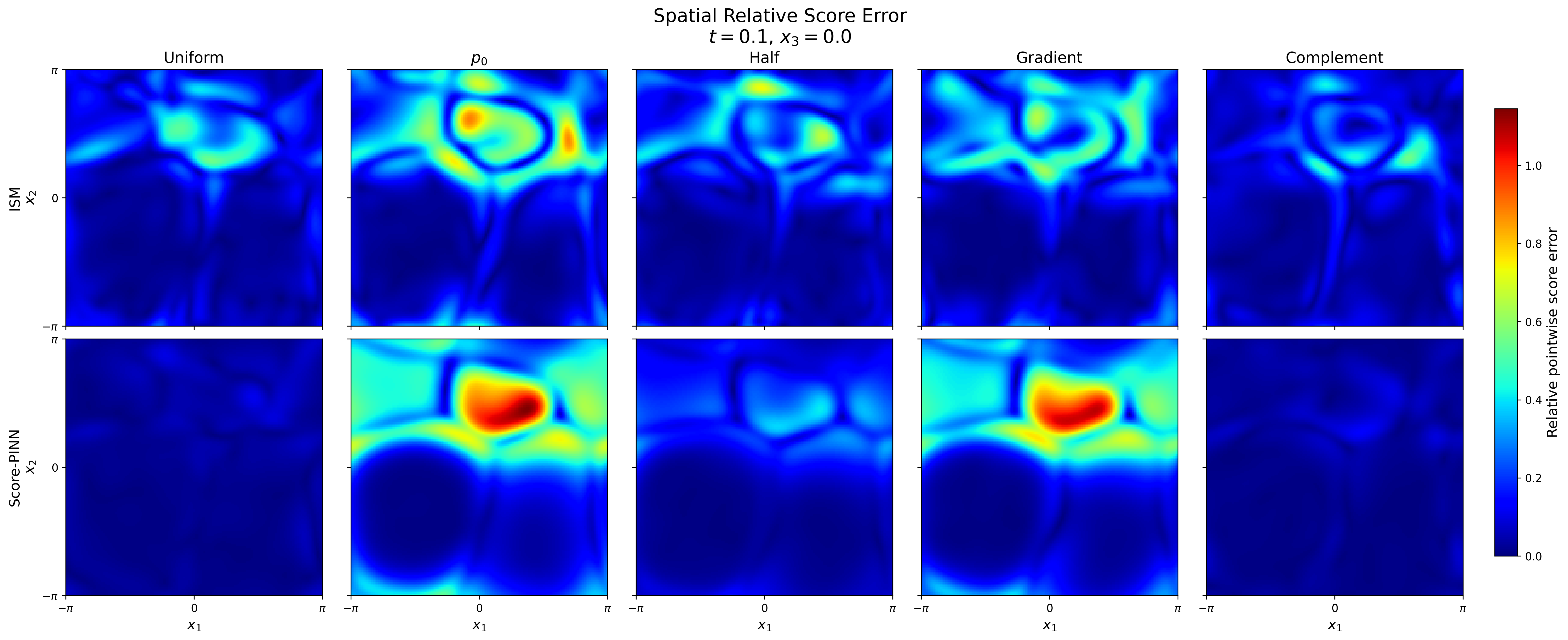}
    \caption{Spatial distributions of the relative pointwise score errors for ISM with BEL supervision (top row) and score-PINN with BEL supervision (bottom row) under five BEL spatial sampling strategies: Uniform, $p_0$, Half, Gradient, and Complement. The errors are evaluated at $t=0.1$ on the spatial slice $x_3=0$, with $N_{\mathrm{BEL}}=1250$. }
    \label{fig:bel_strategies_heatmap_L2_error}
\end{figure}

The relative $L^2$ score errors at $t=0.1$, $0.5$, and $1.0$ are reported in Table~\ref{tab:BEL_strategies_L2_error}. Figure~\ref{fig:bel_strategies_heatmap_L2_error} presents the spatial distributions of the relative pointwise score errors at $t=0.1$ on the spatial slice $(x_1,x_2)$ at $x_3=0$, with the same color scale used across all panels.


\begin{table*}[htp]
    \centering
    \begin{tabular}{llccccc}
        \toprule
        \textbf{Method} & $\boldsymbol{t}$ & \textbf{Uniform} & $\boldsymbol{p_0}$ & \textbf{Half} & \textbf{Gradient} & \textbf{Complement} \\
        \midrule
        ISM & 0.1 & 0.141173 & 0.291683 & 0.190528 & 0.257097 & 0.150052 \\
        ISM & 0.5 & 0.066218 & 0.168876 & 0.105858 & 0.213075 & 0.077449 \\
        ISM & 1.0 & 0.099447 & 0.252584 & 0.139123 & 0.307198 & 0.106942 \\
        \midrule
        Score-PINN & 0.1 & 0.042559 & 0.433131 & 0.137409 & 0.408573 & 0.035162 \\
        Score-PINN & 0.5 & 0.032383 & 0.173649 & 0.043337 & 0.153848 & 0.032058 \\
        Score-PINN & 1.0 & 0.048438 & 0.226847 & 0.064014 & 0.202773 & 0.050974 \\
        \bottomrule
    \end{tabular}
    \caption{Relative spatial $L^2$ score errors under five BEL
    sampling strategies with $N_{\mathrm{BEL}}=1250$, evaluated using
    50,000 uniform spatial points and the respective numerical references.}
    \label{tab:BEL_strategies_L2_error}
\end{table*}

As shown in Figure~\ref{fig:bel_strategies_heatmap_L2_error}, Uniform and Complement sampling produce relatively small errors across most of the spatial domain, whereas $p_0$ and Gradient sampling exhibit more pronounced high-error regions, particularly near the upper projected mode identified in Figure~\ref{fig:marginal_density_GMM4}. Half sampling shows intermediate behavior. These spatial patterns are consistent with the relative $L^2$ errors in Table~\ref{tab:BEL_strategies_L2_error}. For example, at $t=0.1$, the score-PINN error is $0.0426$ under Uniform sampling, compared with $0.4331$ and $0.4086$ under $p_0$ and Gradient sampling, respectively. The corresponding ISM errors are $0.1412$, $0.2917$, and $0.2571$. Overall, these results suggest that adequate spatial coverage, including low-density regions, is important for accurate score estimation over the entire domain.

\subsection{Variance scaling}
\label{app:Variance_scaling}
The Kuramoto process that we chose for the experiment is: 
\begin{equation}
\label{eq:Kuramoto-drift}
     \mathrm{d} \phi_i = 1 + \frac{1}{D} \sum_{j = 1, j \neq i}^{D} \sin(\phi_j - \phi_i) + dW^{i}_t \,.
\end{equation}
It is well known for its synchronizing dynamics \cite{acebron2005kuramoto}. Thus, the Kuramoto drift is given as 
\[
f(x) = \begin{bmatrix}
    1 + \frac{1}{D} \sum_{j = 1, j \neq 1}^{D} \sin(x_j - x_1) + dW^{1}_t\\
    1 + \frac{1}{D} \sum_{j = 1, j \neq 2}^{D} \sin(x_j - x_2) + dW^{2}_t\\
    \dots\\
    1 + \frac{1}{D} \sum_{j = 1, j \neq D}^{D} \sin(x_j - x_D) + dW^{D}_t
\end{bmatrix}
\]
We have
\begin{equation}
     c = -\nabla \cdot f = \frac{1}{D}\sum_{i = 1}^D \sum_{j = 1, j \neq i}^D \cos(\phi_j - \phi_i) \,.
\end{equation}
Thus, the Kuramoto process is defined on torus but not divergence free. Experimenting with the Kuramoto process separates the effect of being divergence-free and being restricted on a torus.

\begin{figure}[H]
    \centering
    \includegraphics[width=\linewidth]{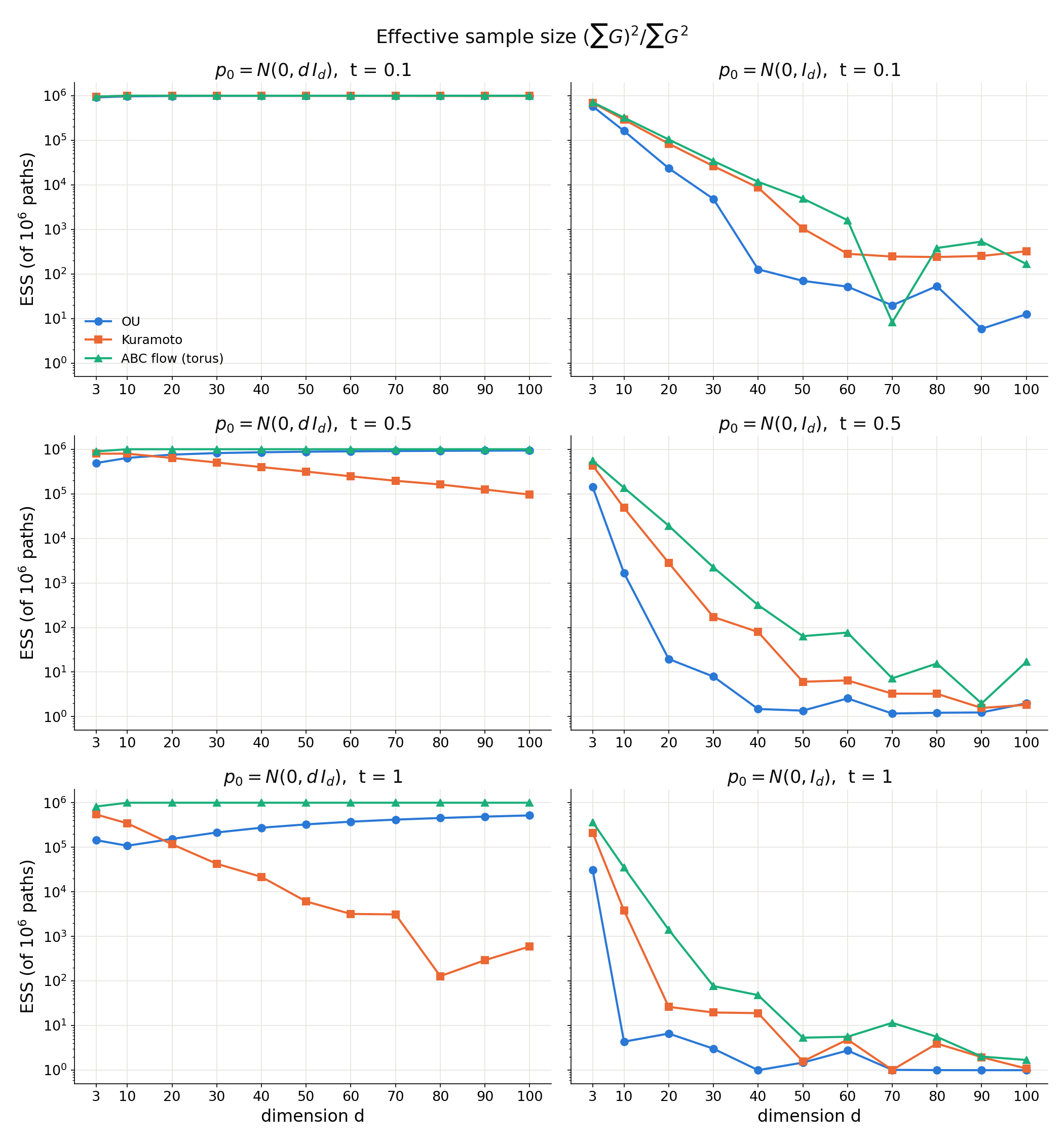}
    \caption{Left column: ESS of $10^6$ paths for the initial distribution $\mathcal{N}(0,dI_d)$. Right column: ESS of $10^6$ paths for the initial distribution $\mathcal{N}(0,I_d)$}
    \label{fig:var-4}
\end{figure}

In addition to the variance check in the main text, we have also computed the effective sample size (ESS) when increasing the sample for Monte Carlo simulation to $10^6$. ESS is given by $= \frac{\sum[G]^2}{\sum[G^2]}$. We have ESS $\approx \frac{N}{R_2}$, where $N$ is the total sample path. Note that the ESS depends $G$ only, so the expansion rate of $\mathbf{H}$ does not affect ESS. The comparison of ESS for three different drifts is compared in Figure \ref{fig:var-4}.

As demonstrated in  Figure~\ref{fig:var-4}, the ESS of the OU process is  decreasing quickly in the case of initial density $\mathcal{N}(0,I_d)$ as either $t$ or $d$ increases. In comparison, the generalized ABC flow has a more robust ESS in higher dimension. The ESS for the Kuramoto process lies in between. 

On the other hand, because the initial density $\mathcal{N}(0, dI_d)$, $p_0$ scales against increasing dimension, the overlap of $\mathcal{N}(0, dI_d)$ with the terminal probability measure is stable with the dimension. This makes the ESS stable for both the OU process and the generalized ABC flow. The exception is Kuramoto process. Because the Kuramoto process is not divergence-free, its $G_t$ contains an exponential term $\exp(\int_0^t c \, \mathrm{d}s)$, making it less stable at larger $t$. The variance scaling experiment confirms the conclusion of our variance analysis in the main text: a divergence-free drift defined on torus has superior variance behavior.
\end{document}